\documentclass[12pt, leqno]{amsart}

\usepackage{
amsmath,
amssymb,
amsthm,
mathrsfs,
amsfonts,
ytableau,
enumitem,
comment,
upgreek,
soul,
yhmath,
mathtools,
shuffle,
kotex
}
\usepackage{tikz,tikz-cd}
\usepackage{latexsym}
\usepackage{mathrsfs}
\usepackage{stackengine}

\usepackage[colorinlistoftodos, textwidth = 2.3cm]{todonotes}

\usepackage[noend]{algpseudocode}
\usepackage[colorlinks=true, pdfstartview=FitV, linkcolor=blue,citecolor=blue, urlcolor=blue]{hyperref}
\makeatletter
\newcommand{\leqnos}{\tagsleft@true\let\veqno\@@leqno}
\newcommand{\reqnos}{\tagsleft@false\let\veqno\@@eqno}
\reqnos
\makeatother

\usepackage{caption}

\ytableausetup{mathmode, boxsize=1.2em}

\newtheorem{theorem}{Theorem}[section]
\newtheorem{proposition}[theorem]{Proposition}
\newtheorem{lemma}[theorem]{Lemma}
\newtheorem{corollary}[theorem]{Corollary}

\newtheorem*{claim*}{Claim}

\theoremstyle{definition}

\newtheorem{example}[theorem]{Example}
\newtheorem{definition}[theorem]{Definition}
\newtheorem{remark}[theorem]{Remark}

\numberwithin{equation}{section} \numberwithin{figure}{section}
\numberwithin{table}{section}

\def\Z{\mathbb Z}

\newcommand{\nc}{\newcommand}

\nc{\La}{\Lambda}
\nc{\ld}{\lambda}
\nc{\rms}{\mathrm{s}}
\nc{\SG}{\mathfrak{S}}
\nc{\SSYT}{\mathsf{SSYT}}
\nc{\stan}{\mathrm{stan}}
\nc{\Span}{\mathrm{span}}
\nc{\comp}{\mathrm{comp}}
\nc{\rmst}{\mathrm{st}}
\nc{\Des}{\mathrm{Des}}
\nc{\set}{\mathrm{set}}
\nc{\wt}{\mathrm{wt}}
\nc{\ch}{\mathrm{ch}}
\nc{\id}{\mathrm{id}}
\nc{\Sym}{\mathrm{Sym}}
\nc{\Qsym}{\mathrm{QSym}}
\nc{\Nsym}{\mathrm{NSym}}
\nc{\sh}{\mathrm{sh}}
\nc{\bfS}{\mathbf{S}}
\nc{\bfm}{\mathbf{m}}
\nc{\hbfS}{\widehat{\mathbf{S}}}
\nc{\bfF}{\mathbf{F}}
\nc{\calS}{\mathcal{S}}
\nc{\hcalS}{\widehat{\mathcal{S}}}
\nc{\alphamax}{\alpha_{\rm max}}
\nc{\brho}{\overline{\rho}}
\nc{\bphi}{\overline{\phi}}
\nc{\calV}{\mathcal{V}}
\nc{\calR}{\mathcal{R}}
\nc{\calG}{\mathcal{G}}
\nc{\tal}{\lambda(\alpha)}
\nc{\tbe}{\widetilde{\beta}}
\nc{\opi}{\overline{\pi}}
\nc{\calP}{\mathcal{P}}
\nc{\rmtop}{\mathrm{top}}
\nc{\rad}{\mathrm{rad}}
\nc{\bfP}{\mathbf{P}}
\nc{\SET}{\mathrm{SET}}
\nc{\SIT}{\mathrm{SIT}}
\nc{\rev}{\mathrm{r}}
\nc{\Th}{\theta}
\nc{\htau}{\widehat{\tau}}
\nc{\mPhi}{\Phi}
\nc{\mphi}{\phi}
\nc{\mPsi}{\Psi}
\nc{\hmPsi}{\widehat{\Psi}}
\nc{\mpsi}{\psi}
\nc{\mGam}{\Gamma}
\nc{\tcd}{\mathtt{cd}}
\nc{\trd}{\mathtt{rd}}
\nc{\trcd}{\mathtt{rcd}}
\nc{\rmr}{\mathrm{r}}
\nc{\rmc}{\mathrm{c}}
\nc{\rmt}{\mathrm{t}}
\nc{\bubact}{\,\scalebox{0.6}{$\bullet$}\,}
\nc{\hbubact}{\,\scalebox{0.6}{$\widehat{\bullet}$}\,}
\nc{\col}{\rm col}
\nc{\row}{\rm row}
\nc{\calE}{\mathcal{E}}
\nc{\calT}{\mathscr{T}}
\nc{\sfT}{\mathsf{T}}
\nc{\calEsa}{\mathcal{E}^\sigma(\alpha)}
\nc{\tauC}{\tau_{\scalebox{0.5}{$C$}}}
\nc{\sytabC}{\sytab_{\scalebox{0.5}{$C$}}}
\nc{\bbfP}{\overline{\bfP}}
\nc{\pr}{\mathsf{pr}}
\nc{\cc}{\mathbf{\mathtt{c}}}
\nc{\Ups}{\Upsilon}
\nc{\pact}{\diamond}
\nc{\tab}{\tau}
\nc{\sytab}{\widehat{\tau}}
\nc{\hatE}{\widehat{E}}
\nc{\hcalE}{\widehat{\calE}}
\nc{\hatC}{\widehat{C}}
\nc{\bal}{{\boldsymbol{\upalpha}}}
\nc{\bbe}{{\boldsymbol{\upbeta}}}
\nc{\bgam}{{\boldsymbol{\upgamma}}}
\nc{\bdel}{{\boldsymbol{\updelta}}}
\nc{\weakcon}{\odot}
\nc{\calB}{\mathcal{B}}
\nc{\ocalB}{\overline{\mathcal{B}}}
\nc{\calM}{\mathcal{M}}
\nc{\calN}{\mathcal{N}}

\newcommand{\la}{\lambda}

\nc{\ldalpha}{\lambda(\alpha)}
\nc{\SRIT}{\mathrm{SRIT}}
\nc{\re}{\mathrm{rev}}
\nc{\otau}{\overline{\tau}}
\nc{\rtop}{{\rm top}}
\nc{\sfc}{\mathsf{c}}
\nc{\calC}{\mathcal{C}}
\nc{\sfr}{\mathsf{r}}
\nc{\tH}{\mathtt{H}}
\nc{\tV}{\mathtt{V}}
\nc{\rpi}{\mathring{\pi}}
\nc{\cpi}{\check{\pi}}
\nc{\frakm}{\mathfrak{m}}
\nc{\fke}{\mathfrak{e}}
\nc{\cact}{\boldsymbol{\mathtt{c}}}
\nc{\bigO}{\mathcal{O}}
\nc{\trms}{\widetilde{\rms}}
\nc{\cont}{\mathrm{cont}}
\nc{\lcm}{\mathrm{lcm}}
\nc{\sfP}{\mathsf{P}}
\nc{\sfhP}{\widehat{\mathsf{P}}}
\nc{\sfM}{\mathsf{M}}
\nc{\bfx}{\mathbf{\mathtt{x}}}
\nc{\bfy}{\mathbf{\mathtt{y}}}

\nc{\bfSsaC}{{\bfS}^\sigma_{\alpha,C}}
\nc{\hbfSsa}{\widehat{\bfS}_\alpha^\sigma}
\nc{\upineq}{\rotatebox{90}{$<$}}
\nc{\downineq}{\rotatebox{270}{$<$}}
\nc{\diagineq}{\rotatebox{135}{$<$}}
\nc{\frakB}{\mathfrak{B}}
\nc{\ttH}{\mathtt{H}}
\nc{\ttV}{\mathtt{V}}
\nc{\rmw}{\mathrm{w}}
\nc{\upsig}{{\boldsymbol{\upsigma}}}
\nc{\bfSsaE}{{\bfS}^\upsig_{\alpha,E}}
\nc{\tst}{\mathtt{st}}

\nc{\ra}{\rightarrow}

\nc{\DIRT}{\mathrm{DIRT}}
\nc{\hpi}{\pi}
\nc{\frakI}{\mathfrak{I}}
\nc{\hfrakI}{\widehat{\mathfrak{I}}}
\nc{\orho}{\overline{\rho}}
\nc{\autotheta}{\uptheta}
\nc{\autopsi}{\uppsi}
\nc{\autophi}{\upphi}
\nc{\autochi}{\upchi}
\nc{\hIM}{\widehat{\frakB}}
\nc{\bfpi}{\boldsymbol{\uppi}}
\nc{\bfopi}{\overline{\boldsymbol{\uppi}}}
\nc{\ofrakB}{\overline{\frakB}}
\nc{\calA}{\mathcal{A}}
\nc{\orb}{\mathtt{orb}}
\nc{\sgn}{\rm{sgn}}

\nc{\sy}[1]{\todo[size=\tiny,color=magenta!10]{#1 \\ \hfill --- So-Yeon}}
\nc{\SY}[1]{\todo[size=\tiny,inline,color=magenta!10]{#1
		\\ \hfill --- So-Yeon}}

\nc{\yt}[1]{\todo[size=\tiny,color=green!10]{#1 \\ \hfill --- Young-Tak}}
\nc{\YT}[1]{\todo[size=\tiny,inline,color=green!10]{#1
		\\ \hfill --- Young-Tak}}
		
\newenvironment{red}{\relax\color{red}}{\relax}
\newenvironment{blue}{\relax\color{blue}}{\hspace*{.5ex}\relax}

\newenvironment{magenta}{\relax\color{magenta}}{\hspace*{.5ex}\relax}

\nc{\ber}{\begin{red}}
\nc{\er}{\end{red}}
\nc{\beb}{\begin{blue}}
\nc{\eb}{\end{blue}}
\nc{\bema}{\begin{magenta}}
\nc{\ema}{\end{magenta}}

\title[Skew Schur polynomials and cyclic sieving phenomenon]{Skew Schur polynomials and cyclic sieving phenomenon}

\author[S.-Y. Lee]{So-Yeon Lee}
\address{Department of Mathematics, Sogang University, Seoul 04107, Republic of Korea}
\email{sylee0814@sogang.ac.kr}

\author[Y.-T. Oh]{Young-Tak Oh}
\address{Department of Mathematics, Sogang University, Seoul 04107, Republic of Korea}
\email{ytoh@sogang.ac.kr}
\thanks{All authors were supported by the National Research Foundation of Korea (NRF) Grant funded by the Korean Government (NRF-2020R1F1A1A01071055).}

\date{\today}
\begin{document}

\maketitle

\begin{abstract}
Let $k$ and $m$ be positive integers
and $\lambda/\mu$ a skew partition.
We compute the principal specialization of the skew Schur polynomials
$s_{\lambda /\mu}(x_1, \ldots, x_{k})$
modulo $q^m-1$ under suitable conditions.
We interpret the results thus obtained 
from the viewpoint of the cyclic sieving phenomenon on semistandard Young skew tableaux of shape $\lambda/\mu$.
As an application, we deal with evaluations of the principal specialization of the skew Schur polynomials at roots of unity.
\end{abstract}

\section{Introduction}
The cyclic sieving phenomenon was introduced by Reiner, Stanton, and White in~\cite{RSW04} as a generalization of Stembridge's $q=-1$ phenomenon~\cite{St94, St96}. 
Let $X$ be a finite set on which a finite cyclic group $C$ acts and let $X(q)$ be a polynomial in $q$ with nonnegative integer coefficients.
For $d\in \Z_{>0}$, let $\omega_d$ be a $d$th primitive root of unity.
We say that the triple $(X, C, X(q))$ exhibits the {\it cyclic sieving phenomenon} if 
$\# X^c = X(\omega_{{\rm o}(c)})$ for all $c\in C$, 
where $X^c$ is the set of fixed points under the action of $c$ 
and ${\rm o}(c)$ is the order of $c$.

We identify a partition with the corresponding Young diagram. 
Given a Young diagram $\lambda$, let $\SSYT_m(\lambda)$ be the set of semistandard Young tableaux of shape $\lambda$ with entries in $\{1,\ldots,m\}$.  
In~\cite{Rho10}, Rhoades proved representation-theoretically that if $\lambda$ is of rectangular shape, the triple
$$
\left(\SSYT_m(\la), ~\langle {\pr} \rangle, ~q^{-\kappa(\la)} s_\la(1,q, \ldots, q^{m-1}) \right) 
$$ 
exhibits  the cyclic sieving phenomenon, where $\pr$ is {\it jeu-de-taquin} promotion,
$\kappa(\la)=\sum_{i\ge 1}(i-1)\la_i$,   
and $s_\la(1,q, \ldots, q^{m-1})$ is the principal specialization of the {\it Schur polynomial} $s_\la(x_1, \ldots, x_m)$.
Since then, there have been several studies on the following issues:
\begin{enumerate}
    \item When $\ld$ is not of rectangular shape, what conditions guarantee the existence
of an action of a finite cyclic group $C$ on $\SSYT_m(\la)$ 
such that the triple 
\[\left(\SSYT_m(\la), ~C, ~q^{-\kappa(\la)} s_\la(1,q, \ldots, q^{m-1}) \right) \] 
exhibits the cyclic sieving phenomenon? 
\item
In (1), can we describe the cyclic group action explicitly? 
\end{enumerate}
For instance, see~\cite{PSMJ21, MBS14, Sam20, OP19, OP21, PW11}.
In particular, it was shown in~\cite{OP19} that for any Young diagram that satisfies the condition $\gcd(|\la|, m)=1$,
the triple 
\begin{align*}
\left(\SSYT_m(\la), ~\langle {\cc} \rangle, ~q^{-\kappa(\la)} s_\la(1, q, \ldots, q^{m-1}) \right) 
\end{align*}
exhibits the cyclic sieving phenomenon.
Here the action $\cc$ arises naturally from the $U_q(\mathfrak{sl}_m)$-{\it crystal} structure of $\SSYT_m(\la)$. 

In this paper, we consider cyclic sieving phenomena associated with semistandard Young skew tableaux.
Unlike tableaux of normal shape, little is known about cyclic sieving phenomena that occur on semistandard Young  skew tableaux.
It has been conjectured by Alexandersson and Amini~\cite[Conjecture 3.4]{AA19} that for any $k, m \in \Z_{>0}$, 
there exists an action of a cyclic group $C_m$ of order $m$ on $\SSYT_k(\la)$ such that 
the triple 
$$
(\SSYT_k(m\la), ~ C_m , ~ s_{m\la}(1, q,  \ldots, q^{k-1}))
$$
exhibits the cyclic sieving phenomenon. 
Here $m\la$ denotes the stretched Young diagram of $\la$ by $m$.
This conjecture turns out to be true in~\cite{OP21}, 
where a crystal-theoretical generalization of this phenomenon is also provided.
Recently, a skew version of this conjecture has been proposed by Alexandersson-Pfannerer-Rubey-Uhlin~\cite[Conjecture 50]{PSMJ21}, which states that  
there is an action of a cyclic group $C_m$ of order $m$ on the set $\SSYT_k(m\lambda/m\mu)$ such that the triple 
\[
(\SSYT_k(m\lambda / m\mu), ~C_m, ~ s_{m\lambda / m\mu}(1, q, \ldots, q^{k-1}))
\]
exhibits the cyclic sieving phenomenon.

The present paper concerns the conjecture given in~\cite[Conjecture 50]{PSMJ21}.
We begin by computing the principal specialization of 
$s_{m\lambda / m\mu}(x_1, \ldots, x_{k})$
modulo $q^m-1$.
For instance, 
in case where $\lambda = (3, 3, 2, 1), \mu = (2, 1), m = 9$, and $k = 4$, 
we observe that
\begin{align*}
s_{9\lambda/9\mu}(1, q, q^2, q^3)
\equiv 54665112 \frac{q^9 - 1}{q - 1} - 3 \frac{q^9 - 1}{q^3 - 1} + 1 \pmod{q^9 - 1}.
\end{align*}
This congruence tells us that the conjecture is not necessarily true.
However, under the condition that $k$ is divisible by $m$, we prove that the conjecture is true. 
To be precise, our first main result is the following. 

\begin{theorem}\label{thm: main theorem}
Let $k$ and $m$ be positive integers
and $\lambda/\mu$ a skew partition.
If $\lambda_i - \mu_i$ is divisible by $m$ for all $i \ge 1$, then
there exists an action of a cyclic group $C_m$ of order $m$ such that the triple 
\[
(\SSYT_{km}(\lambda/\mu), ~C_m, ~s_{\lambda /\mu}(1,q, \ldots, q^{km -1}) )
\]
exhibits the cyclic sieving phenomenon.
\end{theorem}

We also prove that the conjecture is true whenever $\lambda/\mu$ is a border strip, which is our second main result.

\begin{theorem}\label{thm: main theorem2}
Let $k$ and $m$ be positive integers
and $\lambda/\mu$ a border strip. 
If $\lambda_i - \mu_i$ is divisible by $m$ for all $i \ge 1$, then 
there exists an action of a cyclic group $C_m$ of order $m$ such that the triple 
\[
(\SSYT_{k}(\lambda/\mu), ~C_m, ~ s_{\lambda /\mu}(1,q, \ldots, q^{k -1}) )
\]
exhibits the cyclic sieving phenomenon.
\end{theorem}

As an application, we deal with evaluations of the principal specialization of the skew Schur polynomials at roots of unity.
For each $d \, | \, k$, Reiner, Stanton, and White~\cite[Theorem 4.3]{RSW04} provide a formula for $s_{\lambda}(1,q, \ldots, q^{k -1})$ when $q$ is specialized to a $d$th primitive root of unity, which involves
the notion of {\it $d$-cores} and {\it $d$-quotients} of $\lambda$.  
Using a result obtained in the proof of Theorem~\ref{thm: main theorem}, 
we successfully generalize this formula to skew shapes
(see Theorem~\ref{thm: evaluation at roots of unity}).
In particular, the case where $d=k$ is investigated in detail.

This paper is organized as follows. 
In Section~\ref{sec: Preliminaries}, 
we collect the materials required to develop our arguments.
In Section~\ref{Main Theorem 1} and Section~\ref{sec: Main Theorem 2}, the proofs of Theorem~\ref{thm: main theorem} 
and Theorem~\ref{thm: main theorem2} are provided, respectively.
In Section~\ref{sec: Prin-special of skew Schur polynomials at roots of unity},
a skew version of~\cite[Theorem 4.3]{RSW04} is provided.
The final section is devoted to a few remarks for readers' understanding.

\section{Preliminaries}\label{sec: Preliminaries}
\subsection{Partitions, skew partitions, and border strip tableaux}
A {\it partition} is any sequence 
$\lambda=(\lambda_1,\lambda_2, \lambda_3,\ldots)$
of nonnegative integers in decreasing order  
and containing only finitely many nonzero terms.
The {\it length} $\ell(\lambda)$ of $\lambda$ is defined to be the
number of positive parts of $\lambda$.
Unless otherwise specified, we usually omit the zero parts of $\lambda$ and simply write as
$\lambda=(\lambda_1,\lambda_2, \ldots, \lambda_{\ell(\lambda)})$.
Given any positive integer $m$, we set $m\lambda$ to be the stretched partition $(m\lambda_1, m\lambda_2, \ldots, m\lambda_{\ell(\lambda)})$. 
A {\it skew partition} is a pair of partitions $(\lambda, \mu)$ such that the Young diagram of $\lambda$ contains the Young diagram of $\mu$; it is denoted by $\lambda/\mu$. If $\lambda = (\lambda_1, \lambda_2, \ldots)$ and $\mu=(\mu_1, \mu_2, \ldots)$, then the containment of diagrams means that $\lambda_i \geq  \mu_i$ for all $i$.

A {\it border strip} is a connected skew shape with no $2\times 2$ square. Define the {\it height} ${\rm ht(B)}$ of a border strip ${\rm B}$ to be one less than its number of rows.
Let $\alpha=(\alpha_1, \alpha_2, \ldots)$ be a weak composition of $n$.
Define a {\it border-strip tableau} of shape $\ld/\mu$ (where $|\ld/\mu|=n$) and type $\alpha$ to be an assignment of positive integers to the squares of $\ld/\mu$ such that 
\begin{itemize}
\item every row and column is weakly increasing,
\item the integer $i$ appears $\alpha_i$ times, and 
\item the set of squares occupied by $i$ forms a border strip.
\end{itemize}
Define the {\it height} ${\rm ht(T)}$ of a border-strip tableau $T$
to be 
\[{\rm ht}(T)={\rm ht(B_1)}+{\rm ht(B_2)}+ \cdots + {\rm ht(B_k)},\]
where ${\rm B_1}, \ldots, {\rm B_k}$ are the nonempty border strips appearing in $T$. 
For more details, see~\cite[Section 7.17]{99Stan}.

\subsection{The abacus model for partitions, cores and quotients}
Let us review the abacus model introduced by James and Kerber~\cite[Ch. 2.7]{81JK} 
to encode partitions and their quotients, some of which will be slightly modified for the simplicity of the arguments.

\begin{figure}[ht]
\begin{tikzpicture}[baseline = 0mm]
\def \hhh{2.0cm}
\def \hhhh{0.35cm}
\def \hhhhh{8cm}
\def \vvv{-0.7cm}
\def \vvvv{-0.1cm}
\def \ccc{3}
\def \cccc{3}
\def \aaaa{0.2cm}
\def \rdsl{0.14cm}
\def \rdss{0.05cm}
\node[] at (0,0*\vvv) () {\small {\rm runner $0$}};
\node[] at (0+\hhh,0*\vvv) () {\small{\rm runner $1$}};
\node[] at (0+\hhh*2,0*\vvv) () {\small{\rm runner $2$}};

\draw[-,dotted] (0,0+\vvv*6.5) to (0,0.5*\vvv);
\draw[-,dotted] (0+\hhh*1,0+\vvv*6.5) to (0+\hhh*1, 0.5*\vvv);
\draw[-,dotted] (0+\hhh*2,0+\vvv*6.5) to (0+\hhh*2, 0.5*\vvv);

\node[] at (0+\hhhh,0+\vvvv+\vvv*1) (a11) {\tiny 0};
\fill (\hhh*0,\vvv*1) circle (\rdss);
\draw[fill = white] (0+\hhh,0+\vvv*1) circle (\rdsl);
\node[] at (0+\hhh+\hhhh,0+\vvvv+\vvv*1) (a12) {\tiny 1};

\node[] at (0+\hhh*2+\hhhh,0+\vvvv+\vvv*1) (a11) {\tiny 2};
\fill (\hhh*2,\vvv*1) circle (\rdss);

\node[] at (0+\hhhh,0+\vvvv+\vvv*2) (a11) {\tiny 3};
\fill (\hhh*0,\vvv*2) circle (\rdss);
\node[] at (0+\hhh+\hhhh,0+\vvvv+\vvv*2) (a12) {\tiny 4};
\fill (\hhh*1,\vvv*2) circle (\rdss);

\node[] at (0+\hhh*2+\hhhh,0+\vvvv+\vvv*2) (a11) {\tiny 5};
\draw[fill = white] (2*\hhh,0+\vvv*2) circle (\rdsl);

\node[] at (0+\hhhh,0+\vvvv+\vvv*3) (a11) {\tiny 6};
\fill (\hhh*0,\vvv*3) circle (\rdss);
\node[] at (0+\hhh+\hhhh,0+\vvvv+\vvv*3) (a12) {\tiny 7};
\fill (\hhh*1,\vvv*3) circle (\rdss);

\node[] at (0+\hhh*2+\hhhh,0+\vvvv+\vvv*3) (a11) {\tiny 8};
\draw[fill = white] (2*\hhh,0+\vvv*3) circle (\rdsl);

\draw[fill = white] (0,0+\vvv*4) circle (\rdsl);
\node[] at (0+\hhhh,0+\vvvv+\vvv*4) (a11) {\tiny 9};
\draw[fill = white] (0+\hhh,0+\vvv*4) circle (\rdsl);
\node[] at (0+\hhh+\hhhh,0+\vvvv+\vvv*4) (a12) {\tiny 10};

\node[] at (0+\hhh*2+\hhhh,0+\vvvv+\vvv*4) (a11) {\tiny 11};
\fill (\hhh*2,\vvv*4) circle (\rdss);

\node[] at (0+\hhhh,0+\vvvv+\vvv*5) (a11) {\tiny 12};
\fill (\hhh*0,\vvv*5) circle (\rdss);

\node[] at (0+\hhh+\hhhh,0+\vvvv+\vvv*5) (a12) {\tiny 13};
\fill (\hhh*1,\vvv*5) circle (\rdss);

\node[] at (0+\hhh*2+\hhhh,0+\vvvv+\vvv*5) (a11) {\tiny 14};
\draw[fill = white] (0+\hhh*2,0+\vvv*5) circle (\rdsl);

\node[] at (0+\hhhh,0+\vvvv+\vvv*6) (a11) {\tiny 15};
\draw[fill = white] (0,0+\vvv*6) circle (\rdsl);

\fill (\hhh*1,\vvv*6) circle (\rdss);
\node[] at (0+\hhh+\hhhh,0+\vvvv+\vvv*6) (a12) {\tiny 16};
\fill (\hhh*2,\vvv*6) circle (\rdss);
\node[] at (0+\hhh*2+\hhhh,0+\vvvv+\vvv*6) (a11) {\tiny 17};
\end{tikzpicture}
\hspace{6ex}
\begin{tikzpicture}[baseline = 0mm]
\def \hhh{2.0cm}
\def \hhhh{0.3cm}
\def \hhhhh{8cm}
\def \vvv{-0.6cm}
\def \vvvv{-0.1cm}
\def \ccc{3}
\def \cccc{3}
\def \aaaa{0.2cm}
\def \rdsl{0.12cm}
\def \rdss{0.05cm}
\node[] at (0,-\vvv) () {\small {\rm runner $0$}};
\node[] at (0+\hhh,-\vvv) () {\small{\rm runner $1$}};
\node[] at (0+\hhh*2,-\vvv) () {\small{\rm runner $2$}};

\draw[-,dotted] (0,0+\vvv*8.5) to (0,-\vvv/2);
\draw[-,dotted] (0+\hhh*1,0+\vvv*8.5) to (0+\hhh*1,-\vvv/2);
\draw[-,dotted] (0+\hhh*2,0+\vvv*8.5) to (0+\hhh*2,-\vvv/2);

\draw[-, dotted, thick] (-0.3*\hhh, -0.5*\vvv) to (2.5*\hhh, -0.5*\vvv);
\draw[-, dotted, thick] (-0.3*\hhh, 2.5*\vvv) to (2.5*\hhh, 2.5*\vvv);

\node[] at (0+\hhhh,0+\vvvv) (a11) {\tiny 0};
\draw[fill = white] (\hhh*0,\vvv*0) circle (\rdsl);

\node[] at (0+\hhh+\hhhh,0+\vvvv) (a12) {\tiny 1};
\draw[fill = white] (\hhh*1,\vvv*0) circle (\rdsl);

\node[] at (0+\hhh*2+\hhhh,0+\vvvv) (a11) {\tiny 2};
\draw[fill = white] (\hhh*2,\vvv*0) circle (\rdsl);

\node[] at (0+\hhhh,0+\vvvv+\vvv*1) (a11) {\tiny 3};
\fill (0,0+\vvv*1) circle (\rdss);

\node[] at (0+\hhh+\hhhh,0+\vvvv+\vvv*1) (a12) {\tiny 4};
\draw[fill = white] (\hhh*1,\vvv*1) circle (\rdsl);

\node[] at (0+\hhh*2+\hhhh,0+\vvvv+\vvv*1) (a11) {\tiny 5};
\draw[fill = white] (\hhh*2, \vvv*1) circle (\rdsl);

\node[] at (0+\hhhh,0+\vvvv+\vvv*2) (a11) {\tiny 6};
\draw[fill = white] (\hhh*0,\vvv*2) circle (\rdsl);

\node[] at (0+\hhh+\hhhh,0+\vvvv+\vvv*2) (a12) {\tiny 7};
\fill (\hhh*1,\vvv*2) circle (\rdss);

\node[] at (0+\hhh*2+\hhhh,0+\vvvv+\vvv*2) (a11) {\tiny 8};
\draw[fill = white] (\hhh*2,\vvv*2) circle (\rdsl);

\node[] at (0+\hhhh,0+\vvvv+\vvv*3) (a11) {\tiny 0};
\fill (\hhh*0,\vvv*3) circle (\rdss);
\draw[fill = white] (0+\hhh,0+\vvv*3) circle (\rdsl);
\node[] at (0+\hhh+\hhhh,0+\vvvv+\vvv*3) (a12) {\tiny 1};

\node[] at (0+\hhh*2+\hhhh,0+\vvvv+\vvv*3) (a11) {\tiny 2};
\fill (\hhh*2,\vvv*3) circle (\rdss);

\node[] at (0+\hhhh,0+\vvvv+\vvv*4) (a11) {\tiny 3};
\fill (\hhh*0,\vvv*4) circle (\rdss);
\node[] at (0+\hhh+\hhhh,0+\vvvv+\vvv*4) (a12) {\tiny 4};
\fill (\hhh*1,\vvv*4) circle (\rdss);

\node[] at (0+\hhh*2+\hhhh,0+\vvvv+\vvv*4) (a11) {\tiny 5};
\draw[fill = white] (2*\hhh,0+\vvv*4) circle (\rdsl);

\node[] at (0+\hhhh,0+\vvvv+\vvv*5) (a11) {\tiny 6};
\fill (\hhh*0,\vvv*5) circle (\rdss);
\node[] at (0+\hhh+\hhhh,0+\vvvv+\vvv*5) (a12) {\tiny 7};
\fill (\hhh*1,\vvv*5) circle (\rdss);

\node[] at (0+\hhh*2+\hhhh,0+\vvvv+\vvv*5) (a11) {\tiny 8};
\draw[fill = white] (2*\hhh,0+\vvv*5) circle (\rdsl);

\draw[fill = white] (0,0+\vvv*6) circle (\rdsl);
\node[] at (0+\hhhh,0+\vvvv+\vvv*6) (a11) {\tiny 9};
\draw[fill = white] (0+\hhh,0+\vvv*6) circle (\rdsl);
\node[] at (0+\hhh+\hhhh,0+\vvvv+\vvv*6) (a12) {\tiny 10};

\node[] at (0+\hhh*2+\hhhh,0+\vvvv+\vvv*6) (a11) {\tiny 11};
\fill (\hhh*2,\vvv*6) circle (\rdss);

\node[] at (0+\hhhh,0+\vvvv+\vvv*7) (a11) {\tiny 12};
\fill (\hhh*0,\vvv*7) circle (\rdss);

\node[] at (0+\hhh+\hhhh,0+\vvvv+\vvv*7) (a12) {\tiny 13};
\fill (\hhh*1,\vvv*7) circle (\rdss);

\node[] at (0+\hhh*2+\hhhh,0+\vvvv+\vvv*7) (a11) {\tiny 14};
\draw[fill = white] (0+\hhh*2,0+\vvv*7) circle (\rdsl);

\node[] at (0+\hhhh,0+\vvvv+\vvv*8) (a11) {\tiny 15};
\draw[fill = white] (0,0+\vvv*8) circle (\rdsl);

\fill (\hhh*1,\vvv*8) circle (\rdss);
\node[] at (0+\hhh+\hhhh,0+\vvvv+\vvv*8) (a12) {\tiny 16};
\fill (\hhh*2,\vvv*8) circle (\rdss);
\node[] at (0+\hhh*2+\hhhh,0+\vvvv+\vvv*8) (a11) {\tiny 17};
\end{tikzpicture}
\captionsetup{justification=centering}
\caption{The $3$-abaci for $(9^2, 6^3, 4, 1)$ and $(9^2, 6^3, 4, 1)/(2, 1^3)$}
\label{fig: abacus eg}
\end{figure}
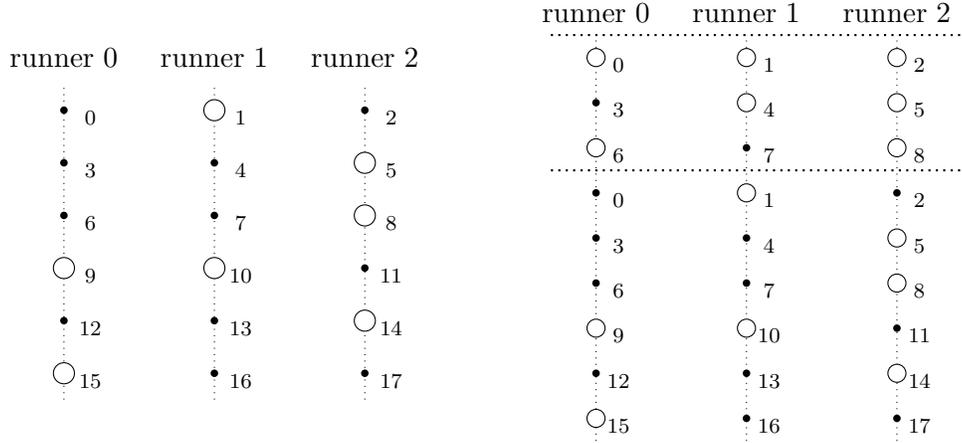

With $d\ge 1$ fixed, we take an abacus with $d$ vertical runners,
numbered $0,1,\ldots, d-1$ from left to right, and we mark positions $0,1,2,\ldots$ on these runners, reading from left to right along successive rows.
Let $\lambda=(\ld_1, \ld_2, \ldots, \ld_\ell)$ be a partition.
For $r \ge \ell$, let $\delta_r$ be the staircase partition $(r-1,r-2,\ldots,1)$
and let $\lambda + \delta_r$ be the sequence
$$
\left(\lambda_1+r-1,\ldots,\lambda_i+(r-i),\ldots,\lambda_{r-1}+1, \ld_r\right),
$$
where $\lambda_i = 0$ if $i > \ell$.
Once $r$ is chosen, the {\it $d$-abacus display for $\ld$}
is obtained by circling the elements of $\lambda + \delta_r$
from the row containing the position labeled by $0$ to the row containing $\lambda_1+r-1$.
It should be mentioned that $r>\ell$ if and only if the position labeled by $0$ is circled.
We call a circled position a {\it bead}
and a uncircled position a {\it non-bead} for simplicity.
See Figure~\ref{fig: abacus eg} for an illustration.

Let $\ld_\emptyset$ be the {\it $d$-core} of $\ld$, which can be obtained by removing
border strips of size $d$ from the diagram of $\ld$ as many as possible.
It is well known that
the abacus display for $\ld_\emptyset$ can be obtained from the abacus display for $\ld$
by moving all the beads as far up their runners as they will go
and then removing the rows containing beads only.
For each $0\le i \le d-1$, let $\lambda^{(i)}$ be the partition whose $j$-part
is the number of non-beads above the $j$th bead from the bottom on runner $i$ of the abacus display for $\ld$.
The sequence $(\lambda^{(0)},\ldots,\lambda^{(d-1)})$
is called {\it the {\it $d$-quotient} of $\lambda$}.

Next, let us introduce the abacus display for skew partitions. 
Let $\lambda/\mu$ be a skew partition with $\ell(\lambda)=\ell$.  
Once $d$ and $r\ge \ell$ are chosen, the $d$-abacus display for $\lambda/\mu$ can be obtained by drawing the $d$-abacus for $\mu$ above the $d$-abacus for $\lambda$.
It should be mentioned that we are using the same $r$ in obtaining the $d$-abaci for $\ld$ and $\mu$. 
In the present paper, $r$ will always be set to be $\ell$.
If the Young diagram of $\lambda^{(i)}$ contains the Young diagram of $\mu^{(i)}$ for all $0\le i \le d-1$, the {\it $d$-quotient} of $\lambda/\mu$ is defined to be the sequence   
$(\lambda^{(0)}/\mu^{(0)},\ldots,\lambda^{(d-1)}/\mu^{(d-1)})$.
Otherwise, we say that the $d$-quotient of $\lambda/\mu$ does not exist.
For more information on the abacus display for skew partitions, we refer the readers to~\cite[Section 5]{PSMJ21}.
\begin{example}
Let $\lambda/\mu=(9^2, 6^3, 4,1)/(2,1^3)$. 
The $3$-abacus for $\lambda/\mu$ can be found in Figure~\ref{fig: abacus eg}.
And, the $3$-quotient of $\lambda/\mu$ is $((4,3)/(1), (2), (2, 1^2))$.
\end{example}

\subsection{Basic results and definitions necessary for proofs}
Here we collect the basic results and definitions necessary to develop our arguments.
The first result to introduce is due to Reiner, Stanton, and White.  
Let $k$, $m$, and $n$ be positive integers.
Suppose that a cyclic group $C_m$ of order $m$ acts {\it nearly freely} on $\{1, \ldots, k\}$,  equivalently it is generated by an element $c \in \SG_k$ whose cycle type is either
\begin{enumerate}
    \item $a$ cycles of size $m$, so that $k = am$ (and $C_m$ acts freely), or
    \item $a$ cycles of size $m$ and one singleton cycle, so that $k = am + 1$
\end{enumerate}
for some positive integer $a$.
Let $X$ be an $n$-multiset of $\{1, \ldots, k\}$, and
\[
X(q) \coloneqq
\begin{bmatrix}
n+k-1 \\
n
\end{bmatrix}_q.
\]
Under this assumption, they prove in~\cite[Theorem 1.1]{RSW04} that 
the triple $(X, C_m, X(q))$ exhibits the cyclic sieving phenomenon.
Furthermore, letting $A_l(k, n)$ denote the coefficient of $q^l$ in 
\begin{align*}
    \begin{bmatrix}
    n+k-1 \\ n
    \end{bmatrix}_q
    \pmod{q^k - 1}
\end{align*}
with the convention that $A_1(1, n) = 1$, 
they show in~\cite[Corollary 1.3]{RSW04} that 
\begin{align}\label{eq: RSW cor}
A_l(k, n) = \sum_{d \mid n,k,l}A_1\left(\frac{k}{d}, \frac{n}{d}\right).
\end{align}
Let $h_n(x_1, \ldots, x_k)$ be the $n$th homogeneous symmetric polynomial in variables 
$x_1, \ldots, x_k$.
It is easy to see that the principal specialization of $h_n(x_1, \ldots, x_k)$ equals $X(q)$,
that is, 
\[
h_n(1, q, \ldots, q^{k-1}) = 
\begin{bmatrix}
n+k-1 \\
n
\end{bmatrix}_q.\]
Thus the theorem~\cite[Theorem 1.1]{RSW04} tells us that 
if a cyclic group of order $m$ act nearly freely on $\{1, \ldots, k\}$,
then  
\[h_n(1, q, \ldots, q^{k - 1}) \equiv \sum_{d \mid m} a_d\frac{q^m - 1}{q^{\frac{m}{d}} - 1} \pmod{q^m - 1}
\]
for some nonnegative integers $a_d$'s.

The second result to introduce is due to Alexandersson and Amini.  

\begin{lemma} {\rm (Alexandersson and Amini~\cite[Theorem 2.7]{AA19})}
\label{lem: necessary and sufficient condition for csp}
Let $f(q)\in \Z_{\ge 0}[q]$ and suppose $f(\omega_m^j)\in \Z_{\ge 0}$ for each $j=1, \ldots, m$, where $\omega_m$ denotes a primitive $m$th root of unity.
Let $X$ be any set of size $f(1)$. 
Then there exists an action of a cyclic group $C$ of order $m$ on $X$ such that $(X, C, f(q))$
exhibits the cyclic sieving phenomenon if and only if for each $d \, | \, m$,
\[\sum_{j \mid d}\mu(d/j)f(\omega_m^j)\ge 0.\]
Here, $\mu$ is the M\"{o}bius function.
\end{lemma}

In Lemma~\ref{lem: necessary and sufficient condition for csp},
suppose that $f(q)$ satisfies the condition that 
\[f(q) \equiv \sum_{d \mid m} a_d\frac{q^m - 1}{q^{\frac{m}{d}} - 1} \pmod{q^m - 1}
\]
for some integers $a_d$'s.
Then it holds that  
\begin{equation}\label{eq: relation to csp polynomials}
\sum_{j \mid d}\mu(d/j)f(\omega_m^j)=da_d
\end{equation}
(for instance, see the proof of~\cite[Theorem 3.2]{OP21}). 
Hence, if $a_d \ge 0$ for all $d \, | \, n$, then 
$(X, C, f(q))$ exhibits the cyclic sieving phenomenon
for a unknown cyclic group action and  
$a_d$ is the number of $C$-orbits with exactly $d$ elements.

Motivated by this observation, we introduce the following definition.
\begin{definition}\label{def: def for csp polynomials}
Let $f(q)\in \Z[q]$.
\begin{enumerate}[label = {\rm (\alph*)}]
\item
We say that $f(q)$ is a {\it pre-CSP polynomial modulo $m$} if 
\begin{align*}
f(q) \equiv \sum_{d \mid m} a_d \frac{q^m - 1}{q^{\frac{m}{d}} - 1} \pmod{q^m - 1}
\end{align*}
for some integers $a_d$'s.
\item 
We say that $f(q)$ is a {\it CSP polynomial modulo $m$} 
if it is a pre-CSP polynomial modulo $m$ with $a_d\geq 0$ for all divisors $d$ of $m$.
\end{enumerate}
\end{definition}

In Definition~\ref{def: def for csp polynomials}(a), the coefficient $a_d$ is denoted by
\[\left[\frac{q^m - 1}{q^{\frac{m}{d}} - 1}\right]f(q)\]
for each divisor  $d$ of $m$ because it is uniquely determined.

\section{Proof of Theorem \protect\ref{thm: main theorem}} \label{Main Theorem 1}
In this section, we provide the proof of Theorem~\ref{thm: main theorem}.
We begin by introducing two important lemmas.
Given positive integers $m$ and $n$, we denote by $(m,n)$ the greatest common divisor of $m$ and $n$
and $[m,n]$ the least common multiple of $m$ and $n$.

\begin{lemma}\label{lem: CSP on hn}
For any positive integers $k$, $m$ and $n$,  
$h_n(1, q, \ldots, q^{km-1})$ is a CSP polynomial modulo $m$.
Furthermore, for each divisor $d$ of $m$, 
\begin{align*}
\left[\frac{q^m - 1}{q^{\frac{m}{d}} - 1}\right] h_n(1,q, \ldots, q^{km-1}) =
\begin{cases}
\left[\frac{q^d - 1}{q - 1}\right] h_{nd/m}(1,q, \ldots, q^{kd-1}) & \text{if $\frac{m}{d}$ divides $n$,}\\
 0 & \text{otherwise.}
\end{cases}
\end{align*}
\end{lemma}

\begin{proof}
By the theorem~\cite[Theorem 1.1]{RSW04}, we can see that 
$h_n(1, q, \ldots, q^{km - 1})$ is a CSP polynomial modulo $m$ as well as a CSP polynomial modulo $km$.
Let $e$ be any divisor of $km$ and set
\[
b_e \coloneqq \left[\frac{q^{km} - 1}{q^{\frac{km}{e}} - 1}  \right]
h_n(1, q, \ldots, q^{km-1}).
\]
In view of~\eqref{eq: RSW cor}, we have that 
\begin{align}\label{eq: coef when k=1}
b_e = \begin{cases}
A_1\left(e, \frac{ne}{km}\right) & \text{if $\frac{km}{e}$ divides $n$,}\\
0 & \text{otherwise.}
\end{cases}
\end{align}
Fix a divisor $d$ of $m$
and let $\calA_d$ be the set of divisors $e$ of $km$ satisfying that $(\frac{km}{e}, m) = \frac{m}{d}$.
Setting 
\[a_d \coloneqq \left[\frac{q^{m} - 1}{q^{\frac{m}{d}} - 1} \right]  h_n(1, q, \ldots, q^{km-1})\]
gives rise to the following relation between $a_d$ and $b_e$:
\begin{align}\label{eq: coef of suit form}
a_d = \sum_{e \in \calA_d} \frac{e}{d}b_e.
\end{align}
In case where $\frac{m}{d}$ divides $n$, 
applying~\eqref{eq: coef when k=1} to this equality yields that 
\begin{align*}
a_d
= \sum_{\substack{e \in \calA_d \\ \frac{km}{e} \mid n}} \frac{e}{d}b_e
+ \sum_{\substack{e \in \calA_d \\ \frac{km}{e} \nmid n}} \frac{e}{d}b_e 
= \sum_{\substack{e \in \calA_d \\ \frac{km}{e} \mid n}} \frac{e}{d}b_e 
= \sum_{\substack{e \in \calA_d \\ \frac{km}{e} \mid n}} \frac{e}{d}A_1\left(e, \frac{ne}{km}\right).
\end{align*}
Otherwise, there is no element $e\in \calA_d$ with $\frac{km}{e}|n$,
and thus $a_d = 0$ by~\eqref{eq: coef when k=1} together with~\eqref{eq: coef of suit form}.

For the second assertion, observe that 
$\calA_d$ is identical to the set of divisors $e$ of $kd$ satisfying that $(\frac{kd}{e}, d) = 1$.
Set
\[
a_d' \coloneqq  \left[ \frac{q^d - 1}{q - 1} \right] h_{nd/m}(1, q, \ldots, q^{kd-1}).
\]
Then we have 
\begin{align*}
a'_d 
=  \sum_{\substack{e \in \calA_d \\ \frac{kd}{e} \mid  \frac{nd}{m} }} \frac{e}{d} A_1\left(e, \frac{ne}{km}\right).
\end{align*}
Now our assertion follows from the fact that, 
for each $e \in \calA_d$, $\frac{km}{e}$ divides $n$ if and only if $\frac{kd}{e}$ divides
$\frac{nd}{m}$.
\end{proof}

\begin{lemma}\label{lem: mutiple of suitable forms}
Let $m$ be a positive integer and let $a, b$ be divisors of $m$. 
Then, 
\begin{equation}\label{eq: modulo multication rule}
\frac{q^m - 1}{q^a - 1} \cdot \frac{q^m - 1}{q^b - 1} \equiv \frac{m}{[a,b]} \cdot \frac{q^m - 1}{q^{(a,b)} - 1} \pmod{q^m - 1}.
\end{equation}
\end{lemma}

\begin{proof}
For simplicity, let $d\coloneqq(a,b)$ and $l\coloneqq[a,b]$.
Then we can write $a = da'$ and $b = db'$ with $(a',b')=1$.
The left hand side of~\eqref{eq: modulo multication rule} equals 
\[\sum_{\substack{0 \leq i < m/a \\ 0 \leq j < m/b}}q^{ai + bj}\]
and the right hand side of~\eqref{eq: modulo multication rule} equals 
\[\frac{m}{l} \sum_{0 \leq k < m/d} q^{dk}.\]
Hence the assertion can be proven by showing that, for each $0 \leq k < \frac md$, the number of solutions to the linear congruence
\begin{equation}\label{linear congruence in consideration}
ai + bj \equiv dk \pmod{m}\quad \text{ with $0 \leq i < \frac ma$  and $0 \leq j < \frac mb$} 
\end{equation}
is given by $\frac ml$. 
To do this, we first note that
the linear congruence
$ai + bj \equiv dk \pmod{m}$ is equivalent to $a'i + b'j \equiv k \pmod{\frac md}$.
Then we consider the set 
\begin{equation}\label{linear congruence in consideration2}
\left\{a'i+b'j \;:\; 0 \le i < \frac ma, \;\; 0 \le j < \frac mb \right\}. 
\end{equation}
For each $0 \leq c < \frac{m}{ab'}$, we claim that the subset  
\[
\left\{a'(i+cb')+b'j \;:\; 0 \le i < b', \;\; 0 \le j < \frac mb \right\}
\]
is a complete system of residues modulo $\frac{m}{d}$.
Note that the number of elements in this set is $\frac md$.   
Suppose that 
\[
a'(i_1 + cb') + b'j_1 \equiv a'(i_2 + cb') + b'j_2 \pmod{\frac{m}{d}},
\]
where $0 \leq i_1, i_2 < b'$ and $0 \leq j_1, j_2 < \frac{m}{b}$.
Then $a'(i_1- i_2) + b'(j_1 - j_2) = \frac{m}{d}k$ for some $k$.
Applying $b' \; | \; \frac{m}{d}$ and $(a',b')=1$ to this equality yields that $b' \; | \; (i_1 - i_2)$. 
Since $-b'< i_1-i_2< b'$, this division relation implies that $i_1=i_2$, thus $j_1 - j_2 = \frac{m}{b}k$. 
It says that $j_1=j_2$ since $-\frac mb < j_1-j_2< \frac mb$, which completes the verification of the claim.

Applying the above claim to \eqref{linear congruence in consideration2}
shows that the number of solutions to \eqref{linear congruence in consideration}
is given by 
\[\frac{\frac{m}{a} \cdot \frac{m}{b}}{\frac{m}{d}} = \frac{m}{l},
\]
as required.
\end{proof}

\begin{remark}\label{rem: sym is suit}
Let $m$ and $k$ be positive integers.
For every partition $\lambda$ with at most $km$ rows,
by applying Lemma~\ref{lem: CSP on hn} and Lemma~\ref{lem: mutiple of suitable forms} repeatedly, we deduce that
$h_\lambda(1, q, \ldots, q^{km - 1})$ 
is a CSP polynomial modulo $m$.
Going further, we deduce that 
the principal specialization of 
every symmetric polynomial in $\Z[x_1,  \ldots, x_{km}]$ is
a pre-CSP polynomial modulo $m$ 
since 
\[\{h_\lambda(x_1, \ldots, x_{km}) : \text{ $\lambda$ a partition with at most $km$ rows}\}\] 
is a basis over $\Z$ of the symmetric polynomials in $x_1, \ldots, x_{km}$.
\end{remark}

To each skew partition $\lambda/\mu$ with  $\ell(\lambda)=l$,
we associate the $l \times l$ matrix $\sfM(\lambda/\mu)$ defined by
\begin{equation}\label{eq: matrix of skew partition} 
\sfM(\lambda/\mu)_{i,j} = \lambda_i - \mu_j + j - i.
\end{equation}

We are now ready to prove Theorem \protect\ref{thm: main theorem}.

\medskip 
{\it Proof of Theorem \protect\ref{thm: main theorem}.}
Let $k$ be an arbitrary positive integer and be fixed throughout the proof.
By Remark~\ref{rem: sym is suit}, we have
\[
s_{\lambda / \mu}(1, q, \ldots, q^{km-1}) \equiv \sum_{d \mid m} c_d\frac{q^m - 1}{q^{\frac{m}{d}} - 1} \pmod{q^m - 1},
\]
where $c_d$'s are integers.
Due to Lemma~\ref{lem: necessary and sufficient condition for csp} together with~\eqref{eq: relation to csp polynomials}, for the assertion, we have only to show  that $s_{\lambda / \mu}(1, q, \ldots, q^{km-1})$ is a CSP polynomial modulo $m$.
This will be accomplished by applying mathematical induction on $m$.
When $m = 1$, there is nothing to prove.  
Suppose that $m \geq 2$ and our assertion holds for all positive integers less than $m$.
For any $n \in \Z$ and any divisor $d$ of $m$,
let 
\begin{equation}\label{eq: congruence form of h}
    a_d(n,m) \coloneqq \left[\frac{q^m - 1}{q^{\frac{m}{d}} - 1} \right] h_n(1, q, \ldots, q^{km - 1}).
\end{equation}
If $n$ is a positive integer, by Lemma~\ref{lem: CSP on hn}, we have
\begin{align}\label{eq: rel b2w coef}
    a_d(n,m) = \begin{cases}
    a_d(\frac{nd}{m}, d) > 0 & \text{if $\frac{m}{d}\, | \,(m,n)$,}\\
    0 & \text{otherwise.}
    \end{cases}
\end{align}
One can also see that $a_d(0,m)=\delta_{d, 1}$ and that when $n$ is a negative integer, $a_d(n,m) = 0$ for all $d \, | \, m$.
Here $\delta_{d, 1}$ denotes the Kronecker delta.

Since every coefficient of $s_{\lambda/\mu}(1, q, \ldots, q^{km - 1})$ is nonnegative, one can see that
$c_m(m) \geq 0$. 
Let $d$ be a divisor of $m$ smaller than $m$, which will be fixed until the end of the proof.
In what follows, we will show that $c_d$ is nonnegative.

Let $\sfM\coloneqq\sfM(\lambda/\mu)$ and $l\coloneqq\ell(\lambda)$.
Consider the $\frac{m}{d}$-abacus display for $\lambda/\mu$ and 
let $(\lambda^{(0)}/\mu^{(0)}, \ldots, \lambda^{(\frac{m}{d} - 1)}/\mu^{(\frac{m}{d} - 1)})$ be the $\frac{m}{d}$-quotient of $\lambda/\mu$.
For $0 \leq t \leq \frac{m}{d} - 1$, let 
\[
\calC_t \coloneqq \{1 \leq i \leq l : \text{the bead for $\lambda_i + l - i$ is on runner $t$ in the abacus display} \}.
\]
Note that for $1 \leq i \leq l$, the bead for $\mu_i + l - i$ is on runner $t$ if and only if $i \in \calC_t$ since $\lambda_i - \mu_i$ is divisible by $m$.
For simplicity, let $l(t) \coloneqq |\calC_{t}|$ and write 
\[
\calC_{t} = \{ c_{t,1} < \cdots < c_{t,l(t)} \}.
\]
Pick any integer $0 \leq t \leq \frac{m}{d} - 1$, which will be fixed until the end of the proof.
Let us introduce three $l(t) \times l(t)$ matrices
$\sfM(t)$, $\sfP(t)$, and $\sfhP(t)$ whose $(i,j)$ entry is given by 
\begin{align*}
&\sfM(t)_{i,j} = \sfM_{c_{t,i}, c_{t,j}},\\
&\sfP(t)_{i,j} = h_{\sfM(t)_{i,j}}(1, q, \ldots, q^{km-1}), \quad \text{and}\\
&\sfhP(t)_{i,j} = \sum_{e \mid d}a_e\left(\sfM(t)_{i,j}, m\right) \frac{q^m - 1}{q^{\frac{m}{e}} - 1}. 
\end{align*}

For $\sigma \in \SG_l$, let
\begin{align*}
f_{\sigma}(q) \coloneqq {\rm sgn}(\sigma) \prod_{1\leq i \leq l} h_{\lambda_i - \mu_{\sigma(i)} + \sigma(i) - i}(1, q, \ldots, q^{km - 1}).
\end{align*}
Let  $\mathfrak{S} \coloneqq \SG_{\calC_{1}} \times \SG_{\calC_{2}} \times \ldots \times \SG_{\calC_{n}}$
and $\SG_{\calC_{t}}$ is the symmetric group on $\calC_{t}$.
In case where $\sigma \in \SG_{l} \setminus \mathfrak{S}$, by applying Lemma~\ref{lem: CSP on hn} and Lemma~\ref{lem: mutiple of suitable forms} to $f_{\sigma}(q)$, we derive that 
\[
\left[\frac{q^m - 1}{q^{\frac{m}{e}} - 1} \right] f_\sigma(q) = 0 \quad \text{for all $e\, | \, d$.}
\] 
Combining this with Jacobi-Trudi formula, we also derive that
\begin{equation}\label{eq: coefficient for some divisors}
\begin{aligned}
    c_e = \sum_{\sigma \in \SG_l} \left[\frac{q^m - 1}{q^{\frac{m}{e}} - 1} \right] f_\sigma(q) 
    = \sum_{\sigma \in \SG} \left[\frac{q^m - 1}{q^{\frac{m}{e}} - 1} \right] f_\sigma(q) = \left[\frac{q^m - 1}{q^{\frac{m}{e}} - 1} \right] \prod_{0 \leq t \leq \frac{m}{d} - 1} \det(\sfP(t))
\end{aligned}
\end{equation}
for all divisors $e$ of $d$, particularly for $e=d$. 
On the other hand,
applying Lemma~\ref{lem: mutiple of suitable forms} to 
\[
\left[\frac{q^m - 1}{q^{\frac{m}{e}} - 1} \right] \sfP(t)_{i,j} = \left[\frac{q^m - 1}{q^{\frac{m}{e}} - 1} \right]\widehat{\sfP}(t)_{i,j}
\quad \text{for all  $1 \leq i,j \leq l$ and $e \, | \, d$}
\]
gives rise to the identity   
\begin{equation}\label{eq: reduction of our assertion}
\left[\frac{q^m - 1}{q^{\frac{m}{e}} - 1} \right] \det(\sfP(t)) = \left[\frac{q^m - 1}{q^{\frac{m}{e}} - 1} \right] \det(\widehat{\sfP}(t)) \quad \text{for all $e \, | \, d$.}
\end{equation}
So our assertion reduces to showing that the right hand side of~\eqref{eq: reduction of our assertion} is always nonnegative.   

We start this by comparing $\sfM(\lambda^{(t)}/\mu^{(t)})$ and $\sfM(t)$.
Note that for $1 \leq i \leq l(t)$,
\begin{align*}
\lambda^{(t)}_i = \frac{d(\lambda_{c_{t,i}} + l - c_{t,i} - t)}{m} - l(t) + i \quad \text{and} \quad
\mu^{(t)}_i = \frac{d(\mu_{c_{t,i}} + l - c_{t,i} - t)}{m} - l(t) + i.
\end{align*}
Therefore, we have 
$\lambda^{(t)} - \mu^{(t)}_j + j - i = \frac{d \sfM(t)_{i,j}}{m}$  for $1 \leq i, j \leq l(t)$.

Since $\lambda^{(t)}_i - \mu^{(t)}_i$ is divisible by  $d$  for all $1 \leq i \leq l(t)$ and $d < m$, by induction hypothesis, we have
\[
s_{\lambda^{(t)}/\mu^{(t)}}(1, q, \ldots, q^{kd-1}) \equiv \sum_{e\mid d}b_e \frac{q^{d} - 1}{q^{\frac{d}{e}} - 1} \pmod{q^{d} - 1},
\]
where $b_e$'s are nonnegative integers.
By replacing $q$ by $q^{\frac{m}{d}}$, we have
\begin{align*}
    s_{\lambda^{(t)}/\mu^{(t)}}(1, q^{\frac{m}{d}}, \ldots, (q^{\frac{m}{d}})^{kd-1}) \equiv \sum_{e\mid d}b_e\frac{q^{m} - 1}{q^{\frac{m}{e}} - 1} \pmod{q^m - 1}.
\end{align*}
On the other hand, since $\frac{d\sfM(t)_{i,j}}{m} = \lambda^{(t)}_i -\mu^{(t)}_j + j-i$,
\eqref{eq: congruence form of h} yields the following congruence:
\begin{align*}
s_{\lambda^{(t)}/\mu^{(t)}}(1, q, \ldots, q^{kd-1}) 
&\equiv \det\left(\sum_{e\mid d}a_{e}\left( \frac{d\sfM(t)_{i,j}}{m}, d\right)\frac{q^{d} - 1}{q^{\frac{d}{e}} - 1}\right)_{1\leq i,j \leq l(t)}
\pmod{q^{d} - 1}.
\end{align*}
This tells us that 
\begin{equation} \label{eq: skew Schur poly for quotients}
\begin{aligned}
s_{\lambda^{(t)}/\mu^{(t)}}(1, q^{\frac{m}{d}}, \ldots, (q^{\frac{m}{d}})^{kd-1})
&\equiv \det\left(\sum_{e\mid d}a_{e}\left(\frac{d\sfM(t)_{i,j}}{m}, d\right)\frac{q^m - 1}{q^{\frac{m}{e}} - 1}\right)_{1\leq i,j \leq l(t)} \\
&\equiv \det\left(\sum_{e \mid d}a_e(\sfM(t)_{i,j}, m)\frac{q^m - 1}{q^{\frac{m}{e}} - 1}\right)_{1\leq i, j \leq l(t)}\\
&\equiv \det\left(\sfhP(t)_{i,j}\right)_{1 \leq i, j \leq l(t)} \pmod{q^m - 1}.
\end{aligned}
\end{equation}
Here the second congruence follows from the formula~\eqref{eq: rel b2w coef}.
As a consequence, 
\[
\left[\frac{q^m - 1}{q^{\frac{m}{e}} - 1} \right] \det(\widehat{\sfP}(t)) = b_e \geq 0 \; \;
\text{for all $e \, | \, d$.}
\]
This completes the proof.
\qed

\begin{remark}
Let $d$ be a divisor of $m$ and $e$ a divisor of $d$.
In the proof of Theorem~\ref{thm: main theorem}, 
combining~\eqref{eq: coefficient for some divisors},~\eqref{eq: reduction of our assertion},
and~\eqref{eq: skew Schur poly for quotients}, 
one can derive the equality
\begin{align}\label{eq: skew ver RSW}
c_e
=\left[\frac{q^m - 1}{q^{\frac{m}{e}} - 1} \right] \prod_{0 \leq t \leq \frac{m}{d} - 1} s_{\lambda^{(t)}/\mu^{(t)}}(1, q^{\frac{m}{d}}, \ldots, (q^{\frac{m}{d}})^{kd-1}).
\end{align}
The idea of obtaining this equality plays an important role in proving Proposition~\ref{prop: skew RSW}.
\end{remark}

\begin{example}
Let $m = 6$, $k = 1$, and $d = 2$.
We here simply write $h_n(q)$ for $h_n(1, q, \ldots, q^5)$ for all $n \in \Z$.
If $\lambda = (13, 10^3, 6)$ and $\mu = (7, 4^3)$, then
$3$-abacus display for $\lambda/\mu$ is
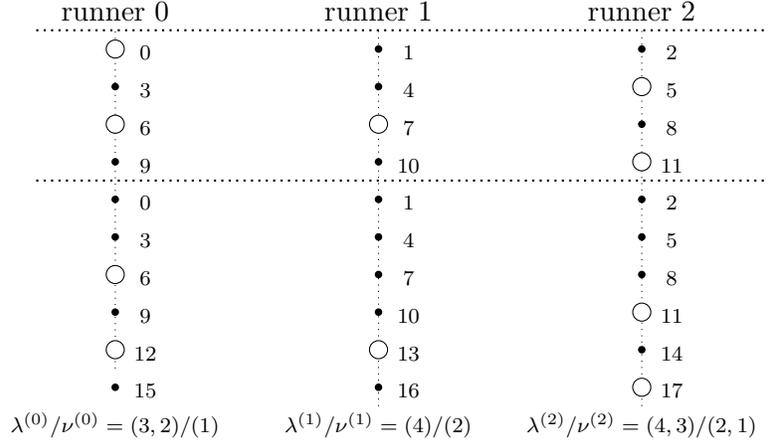
\[
\begin{tikzpicture}[baseline = 0mm]
\def \hhh{3.5cm}
\def \hhhh{0.4cm}
\def \hhhhh{8cm}
\def \vvv{-0.5cm}
\def \vvvv{-0.05cm}
\def \ccc{3}
\def \cccc{3}
\def \aaaa{0.2cm}
\def \rdsl{0.12cm}
\def \rdss{0.05cm}
\node[] at (0,-\vvv) () {\small {\rm runner $0$}};
\node[] at (0+\hhh,-\vvv) () {\small{\rm runner $1$}};
\node[] at (0+\hhh*2,-\vvv) () {\small{\rm runner $2$}};

\draw[-,dotted] (0,0+\vvv*8.5) to (0,-\vvv/2);
\draw[-,dotted] (0+\hhh*1,0+\vvv*9.5) to (0+\hhh*1,-\vvv/2);
\draw[-,dotted] (0+\hhh*2,0+\vvv*9.5) to (0+\hhh*2,-\vvv/2);

\draw[-, dotted, thick] (-0.3*\hhh, -0.5*\vvv) to (2.5*\hhh, -0.5*\vvv);
\draw[-, dotted, thick] (-0.3*\hhh, 3.5*\vvv) to (2.5*\hhh, 3.5*\vvv);

\node[] at (0+\hhhh,0+\vvvv) (a11) {\tiny 0};
\draw[fill = white] (\hhh*0,\vvv*0) circle (\rdsl);

\node[] at (0+\hhh+\hhhh,0+\vvvv) (a12) {\tiny 1};
\fill (\hhh*1,\vvv*0) circle (\rdss);

\node[] at (0+\hhh*2+\hhhh,0+\vvvv) (a11) {\tiny 2};
\fill (\hhh*2,\vvv*0) circle (\rdss);

\node[] at (0+\hhhh,0+\vvvv+\vvv*1) (a11) {\tiny 3};
\fill (0,0+\vvv*1) circle (\rdss);

\node[] at (0+\hhh+\hhhh,0+\vvvv+\vvv*1) (a12) {\tiny 4};
\fill (\hhh*1,\vvv*1) circle (\rdss);

\node[] at (0+\hhh*2+\hhhh,0+\vvvv+\vvv*1) (a11) {\tiny 5};
\draw[fill = white] (\hhh*2, \vvv*1) circle (\rdsl);

\node[] at (0+\hhhh,0+\vvvv+\vvv*2) (a11) {\tiny 6};
\draw[fill = white] (\hhh*0,\vvv*2) circle (\rdsl);

\node[] at (0+\hhh+\hhhh,0+\vvvv+\vvv*2) (a12) {\tiny 7};
\draw[fill = white] (\hhh*1,\vvv*2) circle (\rdsl);

\node[] at (0+\hhh*2+\hhhh,0+\vvvv+\vvv*2) (a11) {\tiny 8};
\fill (\hhh*2,\vvv*2) circle (\rdss);

\node[] at (0+\hhhh,0+\vvvv+\vvv*3) (a11) {\tiny 9};
\fill (\hhh*0,\vvv*3) circle (\rdss);
\fill (0+\hhh,0+\vvv*3) circle (\rdss);
\node[] at (0+\hhh+\hhhh,0+\vvvv+\vvv*3) (a12) {\tiny 10};

\node[] at (0+\hhh*2+\hhhh,0+\vvvv+\vvv*3) (a11) {\tiny 11};
\draw[fill = white] (\hhh*2,\vvv*3) circle (\rdsl);

\node[] at (0+\hhhh,0+\vvvv+\vvv*4) (a11) {\tiny 0};
\fill (\hhh*0,\vvv*4) circle (\rdss);
\node[] at (0+\hhh+\hhhh,0+\vvvv+\vvv*4) (a12) {\tiny 1};
\fill (\hhh*1,\vvv*4) circle (\rdss);

\node[] at (0+\hhh*2+\hhhh,0+\vvvv+\vvv*4) (a11) {\tiny 2};
\fill (2*\hhh,0+\vvv*4) circle (\rdss);

\node[] at (0+\hhhh,0+\vvvv+\vvv*5) (a11) {\tiny 3};
\fill (\hhh*0,\vvv*5) circle (\rdss);
\node[] at (0+\hhh+\hhhh,0+\vvvv+\vvv*5) (a12) {\tiny 4};
\fill (\hhh*1,\vvv*5) circle (\rdss);

\node[] at (0+\hhh*2+\hhhh,0+\vvvv+\vvv*5) (a11) {\tiny 5};
\fill (2*\hhh,0+\vvv*5) circle (\rdss);

\draw[fill = white] (0,0+\vvv*6) circle (\rdsl);

\node[] at (0+\hhhh,0+\vvvv+\vvv*6) (a11) {\tiny 6};
\fill (0+\hhh,0+\vvv*6) circle (\rdss);
\node[] at (0+\hhh+\hhhh,0+\vvvv+\vvv*6) (a12) {\tiny 7};

\node[] at (0+\hhh*2+\hhhh,0+\vvvv+\vvv*6) (a11) {\tiny 8};
\fill (\hhh*2,\vvv*6) circle (\rdss);

\node[] at (0+\hhhh,0+\vvvv+\vvv*7) (a11) {\tiny 9};
\fill (\hhh*0,\vvv*7) circle (\rdss);

\node[] at (0+\hhh+\hhhh,0+\vvvv+\vvv*7) (a12) {\tiny 10};
\fill (\hhh*1,\vvv*7) circle (\rdss);

\node[] at (0+\hhh*2+\hhhh,0+\vvvv+\vvv*7) (a11) {\tiny 11};
\draw[fill = white] (0+\hhh*2,0+\vvv*7) circle (\rdsl);

\node[] at (0+\hhhh,0+\vvvv+\vvv*8) (a11) {\tiny 12};
\draw[fill = white] (0,0+\vvv*8) circle (\rdsl);

\draw[fill = white] (\hhh*1,\vvv*8) circle (\rdsl);
\node[] at (0+\hhh+\hhhh,0+\vvvv+\vvv*8) (a12) {\tiny 13};
\fill (\hhh*2,\vvv*8) circle (\rdss);
\node[] at (0+\hhh*2+\hhhh,0+\vvvv+\vvv*8) (a11) {\tiny 14};

\node[] at (0+\hhhh,0+\vvvv+\vvv*9) (a11) {\tiny 15};
\fill (0,0+\vvv*9) circle (\rdss);

\fill (\hhh*1,\vvv*9) circle (\rdss);
\node[] at (0+\hhh+\hhhh,0+\vvvv+\vvv*9) (a12) {\tiny 16};
\draw[fill = white] (\hhh*2,\vvv*9) circle (\rdsl);
\node[] at (0+\hhh*2+\hhhh,0+\vvvv+\vvv*9) (a11) {\tiny 17};

\node[] at (0, \vvv*10) (a11) {\tiny $\lambda^{(0)}/\nu^{(0)} = (3,2)/(1)$};
\node[] at (0+\hhh*1, \vvv*10) (a11) {\tiny $\lambda^{(1)}/\nu^{(1)} = (4)/(2)$};
\node[] at (0+\hhh*2, \vvv*10) (a11) {\tiny $\lambda^{(2)}/\nu^{(2)} =(4,3)/(2,1)$};
\end{tikzpicture}
\]
and $s_{\lambda /\mu}(1, q, \ldots, q^5) = \det(\sfP)$ with
\begin{align*}
\sfP \coloneqq 
\begin{bmatrix}
{h_{6}(q)} & h_{10}(q) & h_{11}(q) & {h_{12}(q)} & h_{17}(q)\\
h_{2}(q)& {h_{6}(q)} & h_{7}(q) & h_{8}(q) & h_{13}(q)\\
h_{1}(q) & h_{5}(q) & {h_{6}(q)} & h_{7}(q) & h_{12}(q)\\
h_{0}(q) & h_{4}(q) & h_{5}(q) & h_{6}(q) & h_{11}(q)\\
h_{-5}(q) & h_{-1}(q) & h_{0}(q) & h_{1}(q) &
{h_{6}(q)}
\end{bmatrix}.
\end{align*}
Since $\calC_0 = \{3, 5\}$, $\calC_1 = \{2\}$, $\calC_2 = \{1, 4\}$, and $\sfP(t) = \sfP|_{\calC_t \times \calC_t}$, it follows that
\begin{align*}
\sfP(0) = \begin{bmatrix}
h_6(q)& h_{12}(q)\\
h_{0}(q) & h_{6}(q)
\end{bmatrix},
\quad 
\sfP(1) = \begin{bmatrix}
h_6(q)
\end{bmatrix},
\quad \text{and} \quad
\sfP(2) = \begin{bmatrix}
h_6(q)& h_{12}(q)\\
h_{0}(q) & h_{6}(q)
\end{bmatrix}.
\end{align*}
\end{example}

\section{Proof of Theorem~\ref{thm: main theorem2}} \label{sec: Main Theorem 2}
In this section, we provide the proof of Theorem~\ref{thm: main theorem2}.
\begin{lemma}\label{lem: need for ribbon}
Let $k$, $m$, and $n$ be positive integers. 
If $n$ is a multiple of $m$, then
\begin{align*}
    h_n(1,q, \ldots, q^{k-1}) 
    \equiv \sum_{0 \leq j \leq k-1} h_{j}(1, q, \ldots, q^{n-1}) \pmod{q^m - 1}.
\end{align*}
\end{lemma}
\begin{proof}
Note that 
\[ h_n(1,q, \ldots, q^{k-1})=\sum_{T\in \SSYT_k((n))}q^{\omega(T)}, \]
where $\omega(T)=\sum_{1\le i \le k}(i-1) m_i$ and $m_i$ is the number of occurrences of $i$ in $T$.
The decomposition 
\[\SSYT_k((n))=\bigcup_{1\le i \le k}\calA_i\] with 
\[\calA_i= \{ T \in \SSYT_k((n)) : \text{the leftmost entry in $T$ is $i$}\}\]
gives rise to the equality
\begin{align*}
h_n(1,q, \ldots, q^{k-1})
&= \sum_{1 \leq i \leq k} q^{i-1}h_{n-1}(q^{i-1}, \ldots, q^{k-1}). 
\end{align*}
Applying $h_{n-1}(q^{i-1}, \ldots, q^{k-1})=q^{(n-1)(i-1)}h_{n-1}(1, q, \ldots, q^{k-i})$ to the above equality, 
we derive that 
\begin{align*}
    h_n(1,q, \ldots, q^{k-1})
    &\equiv \sum_{1 \leq i \leq k} h_{n-1}(1, q, \ldots, q^{k - i}) \\
    &\equiv \sum_{0 \leq j \leq k-1} h_{n-1}(1, q, \ldots, q^{j})  \quad \text{(by letting $j \coloneqq k-i$)} \\
    &\equiv \sum_{0 \leq j \leq k-1} h_{j}(1, q, \ldots, q^{n-1}) \pmod{q^m - 1}.
\end{align*}
 The third congruence follows from the equality $h_{n-1}(1, q, \ldots, q^{j}) = h_{j}(1, q, \ldots, q^{n-1})$. 
\end{proof}

We are now ready to prove Theorem \protect\ref{thm: main theorem2}.

\medskip 
{\it Proof of Theorem \protect\ref{thm: main theorem2}.}
Let $l \coloneqq \ell(\lambda)$. Since $\lambda/\mu$ is a border strip, we may assume that $\ell(\mu)=l-1$. 
Under this assumption, $\lambda/\mu$ can be naturally identified with the composition 
$\alpha = (\alpha_1, \alpha_2, \ldots, \alpha_l)$ of $|\lambda|-|\mu|$, where  $\alpha_i=\lambda_i-\mu_i$ for $1 \le i \le l$.
From now on, we simply write $r_{\alpha}(1, q, \ldots, q^{k-1})$ and $\sfM(\alpha)$ for 
$s_{\lambda/\mu}(1, q, \ldots, q^{k-1})$ and $\sfM(\lambda/\mu)$, respectively.
By~\eqref{eq: matrix of skew partition}, it holds that 
\begin{align*}
\sfM(\alpha)_{i,j} = \begin{cases}
\alpha_i + \cdots + \alpha_j & \text{if $1 \leq i \leq j \leq l$,}\\
0 & \text{if $2 \leq i \leq l$ and $i - j = 1$,} \\
<0 & \text{otherwise.}
\end{cases}
\end{align*}
Let $\sfP(\alpha; k)$ be the $l \times l$ matrix defined by
\begin{align*}
    \sfP(\alpha; k)_{i,j} = h_{\sfM(\alpha)_{i,j}}(1, q, \ldots, q^{k-1}).
\end{align*}
Combining Jacobi-Trudi formula
with Lemma~\ref{lem: CSP on hn} and Lemma~\ref{lem: mutiple of suitable forms},
we can easily derive that 
\begin{align} \label{eq: csp of ribbon schur}
    r_{\alpha}(1, q, \ldots, q^{k-1}) 
    = \det(\sfP(\alpha; k)) 
    \equiv \sum_{d\mid m}c_d(\alpha, k, m) \frac{q^m - 1}{q^{\frac{m}{d}} - 1} \pmod{q^m - 1},
\end{align}
where $c_d(\alpha, k, m)$'s are integers.
For our purpose, it should be verified that $c_d(\alpha, k, m) \geq 0$ for every divisor $d$ of $m$. 
This will be accomplished by applying induction on $m$.
When $m = 1$, there is nothing to prove.
Suppose that $m \geq 2$ and our assertion holds for all positive integers less than $m$.
Since every coefficient of $ r_{\alpha}(1, q, \ldots, q^{k-1})$ is nonnegative,
it follows that $c_m(\alpha, k, m) \geq 0$.
We now choose an arbitrary divisor $d$ of $m$ smaller than $m$ and 
set $k'\coloneqq\lfloor\frac{d(k-1)}{m}\rfloor + 1$, where $\lfloor\frac{d(k-1)}{m}\rfloor$ denotes the integer part of $\frac{d(k-1)}{m}$.

Let $\widehat{\sfP}(d)$ be the $l \times l$ matrix defined by
\begin{align*}
\widehat{\sfP}(d)_{i,j} = \begin{cases}
 \sum_{0 \leq t \leq k'-1} h_{mt/d}(1, q, \ldots, q^{\sfM(\alpha)_{i,j} - 1}) & \text{if $1 \leq i \leq j \leq l$,}\\
 1 & \text{if $2 \leq i \leq l$,  $j = i - 1$,}\\
 0 & \text{otherwise.}
\end{cases}
\end{align*}
It is clear that  
$\sfP(\alpha;k)_{i,j} = \widehat{\sfP}(d)_{i,j}$ for  $1 \leq j < i \leq l$.
On the other hand, in case where $1 \leq i \leq j \leq l$, Lemma~\ref{lem: CSP on hn} and Lemma~\ref{lem: need for ribbon} imply that  
$\sfP(\alpha;k)_{i,j}$ is a CSP polynomial modulo $m$ satisfying that
\[
\left[ \frac{q^m - 1}{q^{\frac{m}{e}} - 1} \right] \sfP(\alpha; k)_{i,j} = \left[ \frac{q^m - 1}{q^{\frac{m}{e}} - 1} \right] \widehat{\sfP}(d)_{i,j} \quad \text{for all $e \, | \, d$}.
\]
So, from Lemma~\ref{lem: mutiple of suitable forms} it follows that $c_d(\alpha, k, m)$ is equal to $c_d' \coloneqq \left[ \frac{q^m - 1}{q^{\frac{m}{d}} - 1} \right] \widehat{\sfP}(d)$.

In the following, we will see that $c_{d}'\ge 0$.
Let $\beta \coloneqq \left(\frac{d\alpha_1}{m}, \frac{d\alpha_2}{m}, \ldots, \frac{d\alpha_l}{m}\right)$.
For $1 \leq j < i \leq l$, it holds that $\widehat{\sfP}(d)_{i, j} = \sfP(\beta;k')_{i, j}$.
We next consider the case where $1 \leq i \leq j \leq l$. 
Recall that 
\[\sfM(\beta)_{i,j} = \frac{d\sfM(\alpha)_{i,j}}{m}.\]
For each $0 \leq t \leq k'-1$, Lemma~\ref{lem: CSP on hn} shows that
$h_t(1, q^{\frac{m}{d}}, \ldots, (q^{\frac{m}{d}})^{\sfM(\beta)_{i,j} - 1})$ is a CSP polynomial modulo $m$
and 
\[
\left[\frac{q^m - 1}{q^{\frac{m}{e}} - 1} \right] h_t(1, q^{\frac{m}{d}}, \ldots, (q^{\frac{m}{d}})^{\sfM(\beta)_{i,j} - 1})
= 
\left[\frac{q^m - 1}{q^{\frac{m}{e}} - 1} \right] h_{mt/d}(1, q, \ldots, q^{\sfM(\alpha)_{i,j} - 1})
\quad \text{for all $e \, | \, d$.}
\]
Hence $\widehat{\sfP}(d)_{i,j}$ and $\sfP(\beta, k')_{i,j}$
are CSP polynomials modulo $m$ satisfying that
\[
\left[ \frac{q^m - 1}{q^{\frac{m}{e}} - 1} \right] \sfP(\beta; k')_{i,j} = \left[ \frac{q^m - 1}{q^{\frac{m}{e}} - 1} \right] \widehat{\sfP}(d)_{i,j} \quad \text{for all $e \, | \, d$}.
\]
Now, by applying Lemma~\ref{lem: mutiple of suitable forms} to this identity, 
we derive that $c_{d}' = c_d(\beta, k', d)$.
Since $\beta_i$ is divisible by $d$ for all $1 \leq i \leq l$ and $d < m$, 
our induction hypothesis implies that $c_d(\beta, k', d) \geq 0$, as required.
\qed

\begin{remark}
In the proof of Theorem~\protect\ref{thm: main theorem2}, we can additionally observe the following.
\begin{enumerate}[label = {\rm (\alph*)}]
\item
Let $(\lambda^{(0)}/\mu^{(0)}, \ldots, \lambda^{(\frac{m}{d} - 1)}/\mu^{(\frac{m}{d} - 1)})$
be the $\frac{m}{d}$-quotient of $\lambda/\mu$. 
Then $\lambda^{(0)}/\mu^{(0)}$ is the border strip such that $\lambda^{(0)}_i - \mu^{(0)}_i = \beta_i$ for all $1 \leq i \leq l$ and $\lambda^{(t)}/\mu^{(t)} = \emptyset$ for all $1 \leq t \leq \frac{m}{d} - 1$. 

\item 
It holds that $c_d(\alpha, k, m)=c_d(\beta, k', d)$,
which can be viewed as an analogue of Lemma~\ref{lem: CSP on hn}
(for the definition of either side, see~\eqref{eq: csp of ribbon schur}).
\end{enumerate}
\end{remark}

\section{Principal specializations of skew Schur polynomials at roots of unity}
\label{sec: Prin-special of skew Schur polynomials at roots of unity}
Let $\lambda/\mu$ be a skew partition of size $n$.
The skew characters $\chi^{\lambda/\mu}(\nu)$ of the symmetric group $\mathfrak S_n$ 
are then defined implicitly via 
\begin{equation}\label{eq: definition of skew characters}
s_{\ld/\mu}(X)=\sum_{\nu}\chi^{\lambda/\mu}(\nu)\frac{p_{\nu}(X)}{z_{\nu}},
\end{equation}
where the sum is over all partitions $\nu$ of $n$, $z_\nu \coloneqq \prod_j m_j!j^{m_j}$, $m_j$ is the number of parts in $\nu$ equal to $j$, and $p_{\nu}(X)$ is the power-sum symmetric function indexed by  $\nu$.
It can be found in~\cite[Theorem 7.17.3]{99Stan} that  
$$\chi^{\lambda/\mu}(\nu)=\sum_{T}(-1)^{\rm ht(T)},$$
summed over all border-strip tableaux of shape $\ld/\mu$ and type $\nu$. 
In case where $\nu = (d^m)$ with $n=dm$, it was proved in~\cite[Corollary 30]{PSMJ21} that 
$$\chi^{\lambda/\mu}((d^m))=\epsilon|{\rm BST}(\ld/\mu,d)|,$$
where ${\rm BST}(\ld/\mu,d)$ is the set of all border-strip tableaux of shape $\ld/\mu$ and type $(d^m)$
and $\epsilon=(-1)^{\rm ht(T)}$ for any $T\in {\rm BST}(\ld/\mu,d)$.
The following lemma was also presented in the same paper.
\begin{lemma}{\rm (\cite[Lemma 18]{PSMJ21})} \label{lem: condition for existence of quotients}
Let $\ld/\mu$ be a skew shape. 
Then the $d$-quotient of $\ld/\mu$ exists 
if and only if ${\rm BST}(\ld/\mu,d)$ is nonempty.
\end{lemma}

In the following,  
we will show that $\epsilon$ is given by the signature of a particular permutation in  $\mathfrak S_n$, assuming that the $d$-quotient of $\lambda/\mu$ exists.
As above, let $|\lambda/\mu|=dm$ and let $l =\ell(\ld)$.
For $0 \le r \leq d - 1$, let
\begin{align*}
&\ld[r]\coloneqq\{1\le j \le l: \ld_j+(l-j) \equiv r \pmod d\} \; \text{ and }\\
&\mu[r]\coloneqq\{1\le j \le l: \mu_j+(l-j) \equiv r \pmod d\}.
\end{align*}
Since the $d$-cores of $\ld$ and $\mu$ coincide, one sees that 
$\ld[r]$ and $\mu[r]$ have the same number of elements. 
Let us enumerate the elements $\ld[r]$ and $\mu[r]$ in the increasing orders:
\begin{align*}
&\ld[r]=\{ a_{r,1}< a_{r, 2}< \cdots < a_{r, i_r}\},\\
&\mu[r]=\{ b_{r,1}< b_{r, 2}< \cdots < b_{r, i_r}\}.
\end{align*}
Mapping $a_{r,i}$ to $b_{r,i}$ for $1\le i \le i_r$ and for $r=0, 1, \ldots, d-1$, one can obtain a permutation on $\{1, \ldots, l\}$. The resulting permutation is denoted by 
${\bf perm}(\lambda/\mu)$. 
For instance, if $\lambda/\mu=(9^2, 6^3, 4,1)/(2,1^3)$, then 
\begin{align*}
& \ld+\delta_l= (15,14,10,9,8,5,1)\equiv (0,2,1,0,2,2,1) \pmod 3,  \\  
& \mu+\delta_l= (8,6,5,4,2,1,0) \equiv (2,0,2,1,2,1,0) \pmod 3,    
\end{align*}
and thus ${\bf perm}(\lambda/\mu)=2147356$.

Since the $d$-quotient of $\lambda/\mu$ exists, one can obtain $\mu$ from $\ld$ by removing border strips of size $d$.
Suppose that one removes a border strip of size $d$ and height $k$ whose highest and rightmost cell is placed in the $i$th row of $\ld$. 
Let us denote the resulting partition by $\tilde \ld$.
It is easily seen that 
\begin{align*}
\tilde \ld_j
=\begin{cases} 
\ld_j & \text{if $1\le j \le i-1$ or $i+k+1 \le j \le l$,}\\
\ld_{j+1}-1 & \text{if  $i \le j \le i+k-1$,}\\
\ld_{i} -d +k & \text{if $j = i+k$,}
\end{cases}
\end{align*}
and therefore 
\[\tilde \ld +\delta_l=(\ldots, \stackunder{$\underbrace{\ld_{i+1}+l-(i+1)}$}{$i^{\rm th}$},
\ld_{i+2}+l-(i+2), \ldots,\ld_{i+k}+l-(i+k),\stackunder{$\underbrace{\ld_i+(l-i)-d}$}{$(i+k)^{\rm th}$}, \ldots). 
\]
Note that the symmetric group $\mathfrak S_l$ acts on the set of weak compositions of length $l$ by place permutation.
Using this action, one can see that   
\[
\tilde \ld +\delta_l=
\begin{cases}
s_{i+k-1}\cdots s_{i+1}s_i  \cdot (\ld +\delta_l-de_i)& \text{ if $k>0$,}\\
{\rm id}_i \cdot (\ld +\delta_l-de_i)& \text{ if $k = 0$,}\\
\end{cases}
\]
where $e_i \coloneqq (0,\ldots, \stackrel{i{\rm th}}{1}, 0,\ldots, 0)$ and $s_j \coloneqq (j \,\,j+1)$ are simple transpositions of $\mathfrak S_l$. 
In case where $k=0$, the subscript in ${\rm id}_i$ is used to indicate the row from which the border strip of zero height is removed.
One can see that $k$ is equal to the length of $s_{i+k-1}\cdots s_{i+1}s_i$.
Continue this process until we get $\mu+\delta_l$ from $\ld+\delta_l$. 
Composing permutations appearing at each step, one can recover ${\bf perm}(\lambda/\mu)$.
The expressions of ${\bf perm}(\lambda/\mu)$ obtained in this way are in bijection with
the methods of obtaining $\mu$ from $\lambda$ by removing the border strips of size $d$.
Since each method can be uniquely described as a border-strip tableaux of shape $\ld/\mu$ and type $(d^m)$ by numbering the strips in descending order, 
we derive that the set of the expressions of ${\bf perm}(\lambda/\mu)$ obtained in this way is  
one-to-one correspondence with ${\rm BST}(\ld/\mu,d)$.

\begin{example} Let $\lambda/\mu=(9^2, 6^3, 4,1)/(2,1^3)$, $l = 7$, and $d=3$. 
Then $\ld+\delta_l= (15,14,10,9,8,5,1)$
and repeat the above process until we get $\mu + \delta_l$ as follows:
\begin{align*}
& (15,14,10,{\red 9},8,5,1){\stackrel{s_4}{\rightsquigarrow}}  (15,14,10,8,{\red 6},5,1) {\stackrel{s_5}{\rightsquigarrow}}  (15,14,10,8,5,{\red 3},1) {\stackrel{s_6}{\rightsquigarrow}}\\ 
&({\red 15},14,10,8,5,1,0) {\stackrel{s_1}{\rightsquigarrow}} (14,{\red 12},10,8,5,1,0) {\stackrel{s_2}{\rightsquigarrow}} (14,10,{\red 9},8,5,1,0) {\stackrel{s_3}{\rightsquigarrow}} \\ 
& (14,{\red 10},8,6,5,1,0) {\stackrel{s_2}{\rightsquigarrow}}  (14,8,{\red 7},6,5,1,0) {\stackrel{s_4s_3}{\rightsquigarrow}} (14,8,6,{\red 5},4,1,0) {\stackrel{s_4}{\rightsquigarrow}}\\ 
& (14,{\red 8},6,4,2,1,0) {\stackrel{s_2}{\rightsquigarrow}} ({\red 14},6,5,4,2, 1,0){\stackrel{{\rm id}_1}{\rightsquigarrow}}({\red 11},6,5,4,2, 1,0) {\stackrel{{\rm id}_1}{\rightsquigarrow}}\\ 
&(8,6,5,4,2, 1,0)=\mu+\delta_l
\end{align*}
Thus we obtain an expression
\[{\bf perm}(\lambda/\mu)={\rm id}_1 \cdot {\rm id}_1 \cdot s_2 \cdot s_4 \cdot  s_4 s_3 \cdot s_2\cdot s_3 \cdot s_2 \cdot s_1 \cdot s_6 \cdot s_5 \cdot s_4\]
and the corresponding border-strip tableau is 
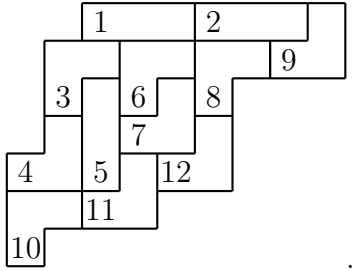
\begin{center}
\def \bp{0.5}
\begin{tikzpicture}
\draw[-, thick] (\bp * 2, \bp * 7) -- (\bp * 9, \bp * 7);
\draw[-, thick] (\bp * 1, \bp * 6) -- (\bp * 2, \bp * 6);
\draw[-, thick] (\bp * 6, \bp * 5) -- (\bp * 9, \bp * 5);
\draw[-, thick] (\bp * 0, \bp * 3) -- (\bp * 1, \bp * 3);
\draw[-, thick] (\bp * 4, \bp * 2) -- (\bp * 6, \bp * 2);
\draw[-, thick] (\bp * 1, \bp * 1) -- (\bp * 4, \bp * 1);
\draw[-, thick] (\bp * 0, \bp * 0) -- (\bp * 1, \bp * 0);
\draw[-, thick] (\bp * 2, \bp * 7) -- (\bp * 2, \bp * 6);
\draw[-, thick] (\bp * 1, \bp * 6) -- (\bp * 1, \bp * 3);
\draw[-, thick] (\bp * 0, \bp * 3) -- (\bp * 0, \bp * 0);
\draw[-, thick] (\bp * 1, \bp * 0) -- (\bp * 1, \bp * 1);
\draw[-, thick] (\bp * 4, \bp * 1) -- (\bp * 4, \bp * 2);
\draw[-, thick] (\bp * 6, \bp * 2) -- (\bp * 6, \bp * 5);
\draw[-, thick] (\bp * 9, \bp * 7) -- (\bp * 9, \bp * 5);

\draw[-, thick] (\bp * 2, \bp * 6) -- (\bp * 8, \bp * 6);

\draw[-, thick] (\bp * 2, \bp * 5) -- (\bp * 3, \bp * 5);
\draw[-, thick] (\bp * 4, \bp * 5) -- (\bp * 5, \bp * 5);

\draw[-, thick] (\bp * 1, \bp * 4) -- (\bp * 2, \bp * 4);
\draw[-, thick] (\bp * 3, \bp * 4) -- (\bp * 4, \bp * 4);
\draw[-, thick] (\bp * 5, \bp * 4) -- (\bp * 6, \bp * 4);

\draw[-, thick] (\bp * 3, \bp * 3) -- (\bp * 5, \bp * 3);
\draw[-, thick] (\bp * 0, \bp * 2) -- (\bp * 3, \bp * 2);

\draw[-, thick] (\bp * 5, \bp * 7) -- (\bp * 5, \bp * 3);

\draw[-, thick] (\bp * 8, \bp * 7) -- (\bp * 8, \bp * 6);
\draw[-, thick] (\bp * 7, \bp * 6) -- (\bp * 7, \bp * 5);

\draw[-, thick] (\bp * 3, \bp * 6) -- (\bp * 3, \bp * 2);

\draw[-, thick] (\bp * 2, \bp * 5) -- (\bp * 2, \bp * 1);

\draw[-, thick] (\bp * 4, \bp * 5) -- (\bp * 4, \bp * 4);

\draw[-, thick] (\bp * 7, \bp * 6) -- (\bp * 7, \bp * 5);

\draw[-, thick] (\bp * 4, \bp * 2) -- (\bp * 4, \bp * 3);

\node at (0.5*\bp, 0.5*\bp) {$10$};
\node at (2.5*\bp, 1.5*\bp) {$11$}; 
\node at (4.5*\bp, 2.5*\bp) {$12$};
\node at (0.5*\bp, 2.5*\bp) {$4$};
\node at (2.5*\bp, 2.5*\bp) {$5$};
\node at (3.5*\bp, 3.5*\bp) {$7$};
\node at (1.5*\bp, 4.5*\bp) {$3$};
\node at (3.5*\bp, 4.5*\bp) {$6$};
\node at (5.5*\bp, 4.5*\bp) {$8$};
\node at (7.5*\bp, 5.5*\bp) {$9$};
\node at (2.5*\bp, 6.5*\bp) {$1$};
\node at (5.5*\bp, 6.5*\bp) {$2$};
\end{tikzpicture}.
\end{center}
\end{example}

\begin{proposition}\label{prop: skew RSW}
Let $d$ be a positive integer. 
If $|\lambda/\mu|=dm$ and the $d$-quotient of $\lambda/\mu$ exists, then 
\[
{\rm sgn}({\bf perm}(\lambda/\mu)) = {\rm sgn} (\chi^{\lambda/\mu}((d^m))).
\]
\end{proposition}

\begin{proof}
Let us obtain an expression of ${\bf perm}(\lambda/\mu)$ in the above way.
The number of simple transpositions appearing in this expression equals the height of the corresponding  
border-strip tableau of shape $\ld/\mu$ and type $(d^m)$.
Now the desired result follows from~\cite[Corollary 30]{PSMJ21}.
\end{proof}

The following theorem is a skew version of~\cite[Theorem 4.3]{RSW04}.
\begin{theorem}\label{thm: evaluation at roots of unity}
Let $d \, | \, N$ and let $\omega_d$ be a primitive $d$th root of unity.
Then $s_{\ld/\mu}(1,\omega_d,\ldots, \omega_d^{N-1})$ is zero unless the $d$-quotient of $\lambda/\mu$ exists, in which case
\[ 
s_{\ld/\mu}(1,\omega_d,\ldots, \omega_d^{N-1})={\rm sgn}(\chi^{\lambda/\mu}((d^m)))\prod_{i=0}^{d-1}s_{\ld^{(i)}/\mu^{(i)}}(\stackunder{$\underbrace{1,1,\ldots, 1}$}{$\frac Nd$}),
\] 
where $|\lambda/\mu|=dm$ and the $d$-quotient of $\ld/\mu$ is $(\ld^{(0)}/\mu^{(0)}, \ldots, \ld^{(d-1)}/\mu^{(d-1)})$.
\end{theorem}

\begin{proof}
In the case where the $d$-quotient of $\lambda/\mu$ does not exist, the proof can be done in the same way as in~\cite[Theorem 4.3]{RSW04}.
Choosing $X = (1, \omega_d, \ldots, \omega_d^{N-1}, 0,0,\ldots)$ 
gives $p_j(X) = 0$ unless $d \, | \, j$.
Hence, in view of~\eqref{eq: definition of skew characters}, one has that 
\begin{equation*}
s_{\ld/\mu}(1,\omega_d,\ldots, \omega_d^{N-1})=\sum_{\nu : d\text{-stretched}}\chi^{\lambda/\mu}(\nu)\frac{p_{\nu}(1,\omega_d,\ldots, \omega_d^{N-1})}{z_{\nu}}.
\end{equation*}
By Lemma~\ref{lem: condition for existence of quotients}, from the condition that the $d$-quotient of $\lambda/\mu$ does not exist it follows that ${\rm BST}(\lambda/\mu, d)$ is empty. 
Thus ${\rm BST}(\lambda/\mu, \nu)$ is also empty.
Combining this with the formula
\[
\chi^{\lambda/\mu}(\nu) = \sum_{{\rm B} \in {\rm BST}(\lambda/\mu, \nu)} (-1)^{\rm ht(B)},
\]
we have $\chi^{\lambda/\mu}(\nu) = 0$ and our assertion follows.

Now suppose that the $d$-quotient of $\lambda/\mu$ exists.
Let $l \coloneqq \ell(\ld)$, $\sigma \coloneqq {\bf perm}(\lambda/\mu)$, and $\sfM \coloneqq \sigma \cdot \sfM(\lambda/\mu)$.
Here $\sigma$ acts on $\sfM(\lambda / \mu)$ by permuting the rows of $\sfM(\lambda / \mu)$.
Then 
\[
s_{\lambda/\mu}(1, q, \ldots, q^{N-1}) = \det(\sfM(\lambda/\mu)) = {\rm sgn}(\sigma) \det(\sfM) = {\rm sgn}(\chi^{\lambda/\mu}((d^m))) \det(\sfM). 
\]
The last equality follows from Proposition~\ref{prop: skew RSW}.
Note that $\det(\sfM)$ is a pre-CSP polynomial modulo $N$.
Write
\[
\det(\sfM) \equiv \sum_{e \mid N} c_e \frac{q^N - 1}{q^{\frac{N}{e}} - 1} \pmod{q^N - 1
 }, 
\]
where $c_e$'s are integers.

Consider the $d$-abacus display for $\lambda/\mu$ and let $(\lambda^{(0)}/\mu^{(0)}, \ldots, \lambda^{(d-1)}/\mu^{(d-1)})$ be the $d$-quotient of $\lambda/\mu$.
For $0 \leq t \leq d-1$, let 
\begin{align*}
\calC_t\coloneqq \{1 \leq i \leq l : \text{the bead for $\mu_i + l - i$ is on runner $t$ of the abacus display}\}.
\end{align*}
By the definition of $\sigma$, one can easily see that the bead for $\lambda_i + l - i$ is on runner $t$ if and only if $\sigma(i) \in \calC_t$.
Therefore, it follows that $\sfM_{i, j}$ is divisible by $d$ if and only if $i, j \in \calC_t$ for some $0 \leq t \leq d-1$.
Using the same argument as in the proof of Theorem~\ref{thm: main theorem}, for each divisor $e \, | \, \frac{N}{d}$, one can derive that
\begin{equation}\label{eq: equalities of some coeffients}
c_e = \left[\frac{q^N - 1}{q^{\frac{N}{e}} - 1}\right] \prod_{0 \leq t \leq d-1} s_{\lambda^{(t)}/\mu^{(t)}}(1, q^d, \ldots, q^{d(\frac{N}{d} - 1)}).
\end{equation}
Using~\eqref{eq: relation to csp polynomials}, one can derive from~\eqref{eq: equalities of some coeffients} that 
\begin{align*}
\det(\sfM)|_{q = \omega_d}
= \sum_{e \mid \frac{N}{d}} e c_e 
&= \prod_{0 \leq t \leq d-1}s_{\lambda^{(t)}/\mu^{(t)}}(1, q^d, \ldots, q^{d(\frac{N}{d} - 1)}) \Bigr|_{q = \omega_d}\\
&= \prod_{0 \leq t \leq d-1}s_{\lambda^{(t)}/\mu^{(t)}}(\stackunder{$\underbrace{1,1,\ldots, 1}$}{$\frac{N}{d}$}).
\end{align*}
\end{proof}

\begin{remark}
(a) Reiner, Stanton, and White proved~\cite[Theorem 4.3]{RSW04} by using the hook-content formula for $s_{\ld}(1,q,\ldots, q^{N-1})$, so one is tempting to
prove Theorem~\ref{thm: evaluation at roots of unity} in a similar manner. 
However, to the best of the authors' knowledge, this formula has not yet extended to general skew shapes (see~\cite[Section 8]{18MPG}).

(b) A similar result has been provided in~\cite[Examples I.5.24]{95Mac}
in a slightly different form.
\end{remark}

The rest of this section is devoted to the case where $d=N$. 
Let $\ell(\ld) \leq N$. 
From~\cite[Theorem 4.3]{RSW04} it follows that 
$s_{\ld}(1,\omega_N,\ldots, \omega_N^{N-1})$ is zero unless the $N$-core of $\lambda$ is empty, in which case
\begin{equation}\label{eq: evaluation when $d=N$}
s_{\ld}(1,\omega_N,\ldots, \omega_N^{N-1})={\rm sgn}(\chi^{\lambda}((N^m)))\quad \text{with  $|\ld|=Nm$.}
\end{equation} 
This evaluation can also be found in~\cite[Examples I.3.17]{95Mac}. 
Recall that the {\it Kostka-Foulkes polynomial} $K_{\ld, \nu}(q)$ is defined via the expansion 
\[s_\ld(x_1,x_2, \ldots, x_N)=\sum_{\nu}K_{\ld, \nu}(q)P_\nu(x_1,x_2,  \ldots, x_N;q),\]
where $P_\nu(x_1,x_2,  \ldots, x_N;q)$ is the Hall-Littlewood polynomial indexed by $\nu$.
And the {\it skew Kostka-Foulkes polynomial} $K_{\ld/\mu, \nu}(q)$ is defined by 
\begin{equation}\label{eq: definition of skew Kostka-Foulkes}
K_{\ld/\mu, \nu}(q)\coloneqq \langle s_{\ld/\mu}(X), ~Q'_{\nu}(X;q) \rangle,
\end{equation}
where $Q'_{\nu}(X;q)=\sum_{\gamma}K_{\gamma, \nu}(q)s_\gamma(X)$ (see~\cite{94DLT}).

Combining~\eqref{eq: evaluation when $d=N$} with~\cite[Theorem 9.14]{94DLT}, for any partition $\ld$ with $|\ld|=Nm$ and $\ell(\ld) \leq N$,
we have 
\begin{equation}\label{eq: evaluation 2 when $d=N$}
K_{\ld, m^N}(\omega_N)=
\begin{cases}
(-1)^{(N-1)m}\,\,{\rm sgn}(\chi^{\lambda}((N^m))) & \text{ if the $N$-core of $\lambda$ is empty,}\\
0 & \text{ otherwise.}
\end{cases}
\end{equation}
We will derive a skew version of this formula. 
To begin with, we observe that  
\[
Q'_{(m^N)}(X;\omega_N)
=(-1)^{(N-1)m} p_N \circ h_m(X)
=(-1)^{(N-1)m}\sum_{\ld \vdash Nm}s_{\ld}(1,\omega_N,\ldots, \omega_N^{N-1})s_\ld(X),
\]
where $p_N \circ h_m$ denotes the plethysm of $h_m$ by $p_N$.
The first equality follows from~\cite[Theorem 9.2]{94DLT} and the second equality follows from~\cite[Examples I.8.7]{95Mac}.
Plugging this identity to~\eqref{eq: definition of skew Kostka-Foulkes} yields that, for any $\ld/\mu$ with $|\ld/\mu|=Nm$,
\begin{equation}\label{eq: evaluation 3 when $d=N$}
K_{\ld/\mu, m^N}(\omega_N)=(-1)^{(N-1)m}s_{\ld/\mu}(1,\omega_N,\ldots, \omega_N^{N-1}).
\end{equation} 

Define a {\it horizontal strip} to be a skew shape with no two squares in the same column.  
The following corollary is the desired skew version of~\eqref{eq: evaluation 2 when $d=N$}.

\begin{corollary} {\rm (cf. \cite[Examples I.5.24]{95Mac})}\label{cor: the case where d=N}
Let $\omega_N$ be a primitive $N$th root of unity.
\begin{enumerate}[label = {\rm (\alph*)}]
\item
$s_{\ld/\mu}(1,\omega_N,\ldots, \omega_N^{N-1})$ is zero unless $\ld^{(i)}/\mu^{(i)}$ is a horizontal strip for all $0\le i\le N-1$, in which case
\[
s_{\ld/\mu}(1,\omega_N,\ldots, \omega_N^{N-1})={\rm sgn}(\chi^{\lambda/\mu}((N^m))) \quad \text{ with } |\ld/\mu|=Nm.
\]
\item For any $\ld/\mu$ with $|\ld/\mu|=Nm$, 
$K_{\ld/\mu, m^N}(\omega_N)$ is zero unless $\ld^{(i)}/\mu^{(i)}$ is a horizontal strip for all $0\le i\le N-1$, in which case
\[K_{\ld/\mu, m^N}(\omega_N)=(-1)^{(N-1)m}\,\,{\rm sgn}(\chi^{\lambda/\mu}((N^m))).
\]
\end{enumerate}
\end{corollary}

\begin{proof}
(a) follows from Theorem~\ref{thm: evaluation at roots of unity} when $d=N$. 

(b) can be obtained by combining (a) with~\eqref{eq: evaluation 3 when $d=N$}.
\end{proof}

\section{Remarks} \label{rem: useful remarks}
We provide a few remarks for readers' understanding.

\begin{enumerate}
\item
As stated in Introduction, the present paper is motivated by the conjecture in~\cite[Conjecture 50]{PSMJ21},
which is equivalent to saying that $s_{m\lambda / m\mu}(1, q, \ldots, q^{k-1})$ is a CSP polynomial modulo $m$. 
It was shown in~\cite{Per19} that this is true in case where $(k, m(|\lambda|-|\mu|))=1$. 
However, one can find a counterexample in case where $k$ is not divisible by $m$ and $(k, m(|\lambda|-|\mu|))\ne 1$. 
For instance, if $\lambda = (3, 3, 2, 1), \mu = (2, 1), m = 9$, and $k = 4$,  using a computer, 
one can easily see that  
\begin{align*}
s_{9\lambda/9\mu}(1, q, q^2, q^{3})\equiv 54665112 \frac{q^9 - 1}{q - 1} - 3 \frac{q^9 - 1}{q^3 - 1} + 1 \pmod{q^9 - 1}.
\end{align*}
On the other hand, in case where $\lambda = (3,3,1)$, $\mu = (2,1)$, $m = 4$, and $k = 6$,
$k$ is not divisible by $m$ but $(k, |m\lambda/m\mu|) = (6,16) = 2 \neq 1$.
In this case, one can see that  
\[
s_{4\lambda/4\mu}(1,q, \ldots, q^{5}) 
\equiv 1576440 \frac{q^4 - 1}{q - 1}
+ 264 \frac{q^4 - 1}{q^2 - 1}
+ 12 \pmod{q^4 - 1}.
\]
It would be very nice to find a sharp condition for the triples $(m\lambda / m\mu, m, k)$ 
such that $s_{m\lambda / m\mu}(1, q, \ldots, q^{k-1})$ is a CSP polynomial modulo $m$. 

\item
Based on~\cite[Theorem 37]{AOE21}, one may naively expect that any stretched skew Schur polynomial $s_{m\lambda/m\mu}(1, q, \ldots, q^{k-1})$ can be a CSP polynomial modulo $m$ when multiplied by an appropriate power of $q$.
However, this immediately turns out to be false.
For instance, as shown in (1), 
$s_{9\lambda/9\mu}(1, q, q^2, q^3)$ is not a CSP polynomial modulo $9$ 
when $\lambda= (3,3,2,1)$, $\mu = (2, 1)$, $m = 9$, and $k = 4$.
In this case, $q^i s_{9\lambda/9\mu}(1, q, q^2, q^3)$ are not CSP polynomials modulo $9$ for all $1 \leq i \leq 8$, and which can be easily checked with the help of a computer. 
Indeed, they are not even pre-CSP polynomials modulo $9$.

\item 
It is well known that there is a one-to-one correspondence 
\[
\SSYT_k(\lambda / \mu) \stackrel{1-1}{\longleftrightarrow} \bigcup_{|\lambda|=|\mu|+|\nu|}\SSYT_k(\nu)^{c^{\lambda}_{\mu \nu}},
\] 
where $c^{\lambda}_{\mu \nu}$ is the Littlewood-Richardson number indexed by $\ld$, $\mu$ and $\nu$ 
(see~\cite[Section 5]{97Ful} or~\cite[Example 3.3]{96GSS}).
However, some $\nu$'s appearing on the right hand side may not be $m$-stretched,
thus \cite[Corollary 3.4]{OP21} cannot be applicable in attacking the conjecture in~\cite[Conjecture 50]{PSMJ21}.
For example, in case where $\lambda = (3,2), \mu = (1)$, $m=4$, and $k = 4$,
there is a one-to-one correspondence between $\SSYT_4(4\lambda/ 4\mu)$ and   
\[
\SSYT_4(11,5) \cup \SSYT_4(12,4) \cup \SSYT_4(8,8) \cup \SSYT_4(10,6) \cup \SSYT_4(9, 7).
\]
\end{enumerate}

\noindent {\bf Acknowledgments.}
The authors thank the anonymous reviewer for providing insightful comments and guidance for further work that resulted in this paper.
They also thank Euiyong Park for helpful discussions.

\bibliographystyle{abbrv}
\bibliography{references}
\end{document}